\pdfoutput=1
\documentclass[final]{IEEEtran}

\usepackage{amsmath,amsfonts,graphicx,amssymb,cite,bm,multirow,multicol,subfigure,mathrsfs}
\usepackage{algorithm,algorithmic}
\usepackage{amscd}
\usepackage[small,nohug,heads=littlevee]{diagrams}
\diagramstyle[labelstyle=\scriptstyle]

\newtheorem{thm}{Theorem}
\newtheorem{lem}{Lemma}

\newtheorem{defn}{Definition}

\newtheorem{exmp}{Example}
\newtheorem{rem}{Remark}

\DeclareGraphicsExtensions{.pdf,.jpeg,.png}

\title{Superposition Frames for Adaptive \\ Time-Frequency Analysis and Fast Reconstruction}

\author{Daniel Rudoy,~\IEEEmembership{Student Member,~IEEE,}
        Prabahan Basu,~\IEEEmembership{Member,~IEEE,}
        and Patrick J. Wolfe,~\IEEEmembership{Senior Member,~IEEE}%
        \thanks{Based upon work supported in part by DARPA Grant~HR0011-07-1-0007, by NIH Grant~P01~CA134294-01, and by NSF Grant~DMS-0652743 and a Graduate Research Fellowship.  Preliminary versions of material in parts of Sections~\ref{sec:Ola} and~\ref{sec:algSpecifics} first appeared in the 2008 IEEE International Conferences on Acoustics, Speech, and Signal Processing (ICASSP)~\cite{Rudoy2008}.
         \newline \indent The authors are with the Statistics and
         Information Sciences Laboratory, Harvard University, Cambridge, MA
         02138 USA (e-mail: \{rudoy, pbasu, patrick\}@seas.harvard.edu).
         }
\vspace{-\baselineskip}%
}

\begin{document}

\maketitle

\begin{abstract}

In this article we introduce a broad family of adaptive, linear time-frequency representations termed superposition frames, and show that they admit desirable fast overlap-add reconstruction properties akin to standard short-time Fourier techniques.  This approach stands in contrast to many adaptive time-frequency representations in the existing literature, which, while more flexible than standard fixed-resolution approaches, typically fail to provide for efficient reconstruction and often lack the regular structure necessary for precise frame-theoretic analysis.  Our main technical contributions come through the development of properties which ensure that our superposition construction provides for a numerically stable, invertible signal representation.  Our primary algorithmic contributions come via the introduction and discussion of specific signal adaptation criteria in deterministic and stochastic settings, based respectively on time-frequency concentration and nonstationarity detection.  We conclude with a short speech enhancement example that serves to highlight potential applications of our approach.
\end{abstract}%


\section{Introduction}
\label{sec:intro}

\IEEEPARstart{O}{vercomplete} short-time Fourier methods are frequently used to analyze the time-varying spectral content of discrete-time waveforms $x[t]$ arising in a variety of signal processing applications. Since the choice of localizing window function effectively controls the balance between time and frequency resolution \emph{a priori}, standard representations cannot modulate this trade-off to adapt to the local spectral content of $x[t]$.  Over the past two decades, this shortcoming has motivated the development of various linear and nonlinear adaptive time-frequency analysis methods~\cite{Jones1990, JonesBaraniuk94, Jones1995, Czerwinski1997, Goodwin1998, Mallat99, Kwok2000, Wolfe2001, Dorfler02, Wolfe04a, Jaillet2007,Djurovic03,Nesbit2009}, in applications ranging from biomedical engineering~\cite{Durka2001} to radar signal analysis~\cite{Rothwell1998} and speech processing~\cite{Heusdens2006}.

Despite the recognized importance of overcomplete signal-adaptive time-frequency analysis, the above methods generally fail to admit \emph{fast reconstruction} of the signal $x[t]$ from its time-frequency representation, by which we mean any non-iterative method that avoids direct (pseudo-)inversion of the corresponding analysis operator.  While approaches such as modulated lapped transforms for audio coding~\cite{malvar1990lte}, wavelet packet decompositions via best basis~\cite{Coifman1992, WesfriedWickerhauser93} and adaptive segmentation via dynamic programming~\cite{RamchandranOrchard97, PrandoniVetterli00, NiamutHeusdens05, HeusdensJensen2005, RodbroHeusdens06}  can lead to flexible tilings of the time-frequency plane, the general goal of efficient reconstruction from signal-adaptive, \emph{overcomplete time-frequency} representations remains an open problem. This issue is particularly important, given the recent interest in oversampled, modulated filter banks~\cite{Cvetkovic1998a, Bolcskei1998, Kovacevic2007a}.

In this article we introduce a broad family of adaptive, linear time-frequency representations that admit a fast overlap-add reconstruction property akin to standard short-time Fourier techniques.  We do so by adapting a given discrete Gabor frame to an observed signal $x[t]$ via superpositions of neighboring translates of a single window function, to yield the superposition frames of the article title.  Related procedures include the multi-window constructions of~\cite{Zibulski1997}, in which multiple systems are defined on the same time-frequency lattice; and the multi-Gabor expansions of~\cite{Li1999}, in which multiple time lattices and windows are employed. However, neither of these schemes treats the use of subset selection to achieve a signal-adaptive system in the manner of the present article.  More recent approaches~\cite{Wolfe2001, Dorfler02, Wolfe04a, Jaillet2007} address subset selection from a Gabor frame or union of Gabor frames, but do not consider the structure of the corresponding canonical dual.

A very recent approach in this direction is the study of general nonstationary Gabor frames~\cite{JailletDorfler09}, and indeed our contribution can be viewed as one possible instantiation of this framework. However, as we show below, the additional structure induced by our superposition construction yields several important properties, including, among other results, a \emph{preservation} of the lower frame bound of the original Gabor frame, a generalized constant overlap-add property that avoids the explicit computation of dual windows, and a means of generating new families of adaptive lapped frames.

The article is organized as follows. We begin by reviewing the short-time Fourier transform and Gabor systems on $\mathbb{C}^L$ in Section~\ref{sec:STFT_Gabor}.  Next, we introduce superposition windows in Section~\ref{sec:superposition} and use them to construct superposition systems in Section~\ref{sec:adaptFramework}. In Section~\ref{sec:adaptiveFrameAnalysis} we prove that the resultant systems are in fact frames for $\mathbb{C}^L$ and study their frame-theoretic properties, and in Section~\ref{sec:synthesis} we establish fast reconstruction via an analysis of the corresponding frame operator. In Section~\ref{sec:algSpecifics}, we give examples of signal-dependent adaptation algorithms and illustrate their application to superposition systems.  We conclude with a brief discussion in Section~\ref{sec:discussion}. %

\section{Preliminaries}
\label{sec:STFT_Gabor}

We first review some well-known properties of Gabor frames~\cite{Qiu1995, christensen2003ifa, Kovacevic2007a} and discuss their relationship to short-time Fourier analysis.  We take as our setting the space $\mathbb{C}^L$, and interpret its members as discrete-time $L$-periodic signals $x \in \ell^2(\mathbb{Z}_L)$, with $\mathbb{Z}_L$ denoting the integers $\mathbb{Z}$ modulo $L$. The \emph{short-time Fourier transform} (STFT) on $\mathbb{C}^L$ uses a well-concentrated window function in order to localize $x$ in time prior to the
analysis of its frequency content.

\begin{defn}[Short-Time Fourier Transform]\label{def:STFT} Fix a window $w \in \mathbb{C}^L$ and time-frequency lattice constants $a, b >0$ that divide $L$, with $a$ an integer, and define $M,N: Na = Mb = L$. Then for the $m$th frequency bin index and $n$th window shift, with $m \in \mathbb{Z}_M$ and $n \in \mathbb{Z}_N$, the Gabor or subsampled short-time Fourier transform $X[m,n]$ of $x \in \mathbb{C}^L$ is given by
\begin{equation}
    \label{eq:STFT}
        X[m,n] \triangleq \sum_{t=0}^{L-1} x[t]\overline{w[t-na]e^{2\pi i mbt/L}}
        \text{,}
\end{equation}
\end{defn}
where $i = \sqrt{-1}$ and $\overline{\,\cdot\,}$ denotes complex conjugation. The expression of~\eqref{eq:STFT} can be viewed as a set of inner products of $x$ with $NM$ time-frequency shifts of the chosen window $w$. To realize this correspondence, and to set notation, we introduce explicit translation and modulation operators as follows.

\begin{defn}[Translation and Modulation Operators]\label{def:transModOprs}
Let the translation and modulation operators $\mathcal{T}$ and $\mathcal{M}$ be defined as maps from $\mathbb{C}^L$ to itself acting according to:
\begin{equation*}
    \mathcal{T}_{na}w[t] \triangleq w[t-na], \quad
    \mathcal{M}_{mb}w[t] \triangleq w[t]\,e^{2\pi i mbt/L} \text{.}
\end{equation*}
\end{defn}

Through the action of these operators, time-frequency shifts of the chosen window $w \in \mathbb{C}^L$ may be indexed as
\begin{equation}
    \label{eq:GaborFrame}
    \phi_{m,n}[t] \triangleq \mathcal{M}_{mb}\mathcal{T}_{na}w[t]
    , \quad m \in \mathbb{Z}_M, n \in \mathbb{Z}_N
    \text{,}
\end{equation}
and one speaks of a Gabor system $\mathscr{G}(w,a,b) = \{\phi_{m,n}\}$. In order to ensure a reconstruction property for any $x$ from its subsampled short-time Fourier transform $X[m,n]$, the Gabor system $\mathscr{G}(w,a,b)$ must form a \emph{frame} for $\mathbb{C}^L$ as follows.

\begin{defn}[Gabor Systems and Frames]\label{def:frameCond}
A denumerable set $\{ \phi_{m,n} \}$ of vectors comprising time-frequency shifts of a single window function $w \in \mathbb{C}^L$ is called a \emph{Gabor system}, and is said to be a \emph{Gabor frame} for $\mathbb{C}^L$ if there exist constants $0 < A \leq B < \infty$ termed \emph{frame bounds} such that:
\begin{equation}
    \label{eq:frameExpansion}
    \textstyle
   \forall x \in  \mathbb{C}^L,\, A\|x\|^2 \le \sum_{m,n} |\langle x,\phi_{m,n} \rangle|^2 \le B \|x\|^2 \text{,}
\end{equation}
with inner product $\langle x,\phi_{m,n} \rangle \triangleq \sum_{t=0}^{L-1} x[t] \overline{\phi_{m,n}[t]} = X[m,n]$.
\end{defn}

An upper frame bound $B$ for~\eqref{eq:frameExpansion} is guaranteed whenever the set $\{ \phi_{m,n} \}$ is finite, and so the existence of a lower frame bound $A > 0$,  for a finite Gabor system $\mathscr{G}(w,a,b)$, is equivalent to the requirement that its elements span $\mathbb{C}^L$. This occurs \emph{if and only if} the \emph{frame operator} is of full rank.
\begin{defn}[Gabor Frame Operator]\label{def:frameOpr}
Let $\mathscr{G}(w,a,b) = \{ \phi_{m,n} \}$ be a Gabor system on $\mathbb{C}^L$, and define the \emph{frame operator} $S: \mathbb{C}^L \rightarrow \mathbb{C}^L$ through its action on $x$ as $Sx = \sum_{m,n} \langle  x, \phi_{m,n} \rangle \phi_{m,n}$.  Then $S$ is represented by the $L \times L$ symmetric and positive semi-definite matrix with entries
\begin{equation}
    \label{eq:Spointwise}
    S[t,t'] \triangleq \sum_{m=0}^{M-1} \sum_{n=0}^{N-1} \mathcal{M}_{mb} \mathcal{T}_{na} w[t] \overline{\mathcal{M}_{mb}\mathcal{T}_{na}w[t']} \text{.}
\end{equation}
\end{defn}

\begin{rem}[Strict Positive-Definiteness of Frame Operator]\label{rem:SposDef}
By Definition~\ref{def:frameOpr}, the frame condition of~\eqref{eq:frameExpansion} is equivalent to \emph{strict} positive definiteness of $S$ and hence a necessary condition is that $MN \geq L$ (i.e., $ab \leq L$). Moreover, the minimal and maximal eigenvalues of $S$ yield optimal frame bounds, since~\eqref{eq:frameExpansion} may be expressed as the requirement that $A \langle x,x \rangle \leq \langle Sx,x \rangle \leq B\langle x,x \rangle, \, \forall x \in  \mathbb{C}^L$.
\end{rem}

The frame condition of~\eqref{eq:frameExpansion} in turn implies the following reconstruction property:
\begin{equation*}
\forall x \in \mathbb{C}^L, t \in \mathbb{Z}_L, \,\, x[t] =
\sum_{m=0}^{M-1} \sum_{n=0}^{N-1} \langle x, \phi_{m,n} \rangle
\widetilde{\phi}_{m,n}[t] \text{,}
\end{equation*}
where the elements $\{ \widetilde{\phi}_{m,n} \}$ comprise a (not necessarily unique) \emph{dual frame}. However, to each frame may be associated a unique \textit{canonical} dual, whose elements are given by the action of the frame operator inverse $S^{-1}$ on each $\phi_{m,n}$. Moreover, in the Gabor setting, this canonical dual takes the form of another Gabor system $\mathscr{G}(\widetilde{w},a,b)$, with $\widetilde{w} \triangleq
S^{-1}w$.

Any $S$ can be written as a sum of outer products of each frame vector with itself; from~\eqref{eq:Spointwise} via the orthogonality
relation
\begin{equation*}
    \sum_{m=0}^{M-1} e^{2\pi i mb(t-t')/L} =
         \begin{cases}
            M & \text{when $M$ divides $t-t'$,} \\
            0 & \text{otherwise} \text{,}
         \end{cases}
\end{equation*}
we obtain the so-called Walnut representation~\cite{HeilWalnut89} of a Gabor frame operator $S$, which will be used repeatedly
throughout.

\begin{defn}[Discrete Walnut Representation]\label{def:Walnut}
Denote by $M \backslash (t-t')$ the condition that $M$ divides $t-t'$, and by $\mathbb{I}_{M \backslash (t-t')}[t-t']$ the corresponding indicator function on $\mathbb{Z}_L$.  Then the frame operator $S$ of a finite Gabor system $\mathscr{G}(w,a,b)$ has banded structure, and satisfies the entrywise relation
\begin{equation}
    \label{eq:Walnut}
    S[t,t'] = \mathbb{I}_{M\backslash(t-t')}[t-t'] \cdot M\!\sum_{n=0}^{N-1}\mathcal{T}_{na}w[t]\overline{\mathcal{T}_{na}w[t']} \text{.}
\end{equation}
\end{defn}

\begin{rem}[Covering Condition]\label{rem:covCond}
Note that if $\mathscr{G}(w,a,b)$ is a frame for $\mathbb{C}^L$, then~\eqref{eq:Walnut} implies that the covering condition
\begin{equation}\label{eq:covCond}
\sum_{n=0}^{N-1} |w[t-na]|^2 > 0, \forall \, t \in \mathbb{Z}_L
\end{equation}
must be satisfied, since a necessary condition for positive definiteness of $S$ is that its diagonal entries are positive.
\end{rem}

\begin{rem}[Window Length as Distinct from Support] \label{rem:windowlength}
The \emph{support} of $w \in \mathbb{C}^L$ refers to the set of indices $t$ for which $w[t] \neq 0$, with $|\operatorname{supp}(w)|$
its cardinality.  Bearing in mind the summands of~\eqref{eq:Walnut} and~\eqref{eq:covCond}, define the \emph{length} of $w$ by
\begin{equation}\label{eq:winLen}
    \operatorname{len}(w) \triangleq |\operatorname{supp}(w)|
\end{equation}
if $\operatorname{supp}(w)$ is contiguous, as is often the case in practice, and $\min_{n\in\mathbb{Z}_L} (\max_{t,t' \in \mathbb{Z}_L}
|t-t'+1| \!:\! \mathcal{T}_n(w[t]\overline{w[t']}) \neq 0 )$
otherwise.
\end{rem}

\begin{rem}[Diagonal Frame Operator]\label{rem:diagFrameOpr}
It follows from~\eqref{eq:Walnut} and~\eqref{eq:winLen} that $S$ is diagonal if $M \geq \operatorname{len}(w)$, since $\mathcal{T}_{na}w[t]\overline{\mathcal{T}_{na}w[t']} = 0$ for all $|t-t'| \geq \operatorname{len}(w)$, including those for which $M$ divides $t-t'$. In turn, this implies efficient computation of the dual frame $\{\mathcal{M}_{mb}\mathcal{T}_{na}\widetilde{w}\}$, with $\widetilde{w} = S^{-1}w$ obtained via element-wise division of $w[t]$ by $S[t,t] = M\sum_{n=0}^{N-1} |w[t-na]|^2$.  In this case the condition of~\eqref{eq:covCond} is sufficient to guarantee the frame condition of~\eqref{eq:frameExpansion}.
\end{rem}

\begin{figure*}[t]
  \centering
\subfigure[\label{fig:windowExamples}Superposition windows $w_r$, for $r = 0, 1, 2, 4$ merges of neighboring translates.]{\hfill
  \includegraphics[width=.29\textwidth]{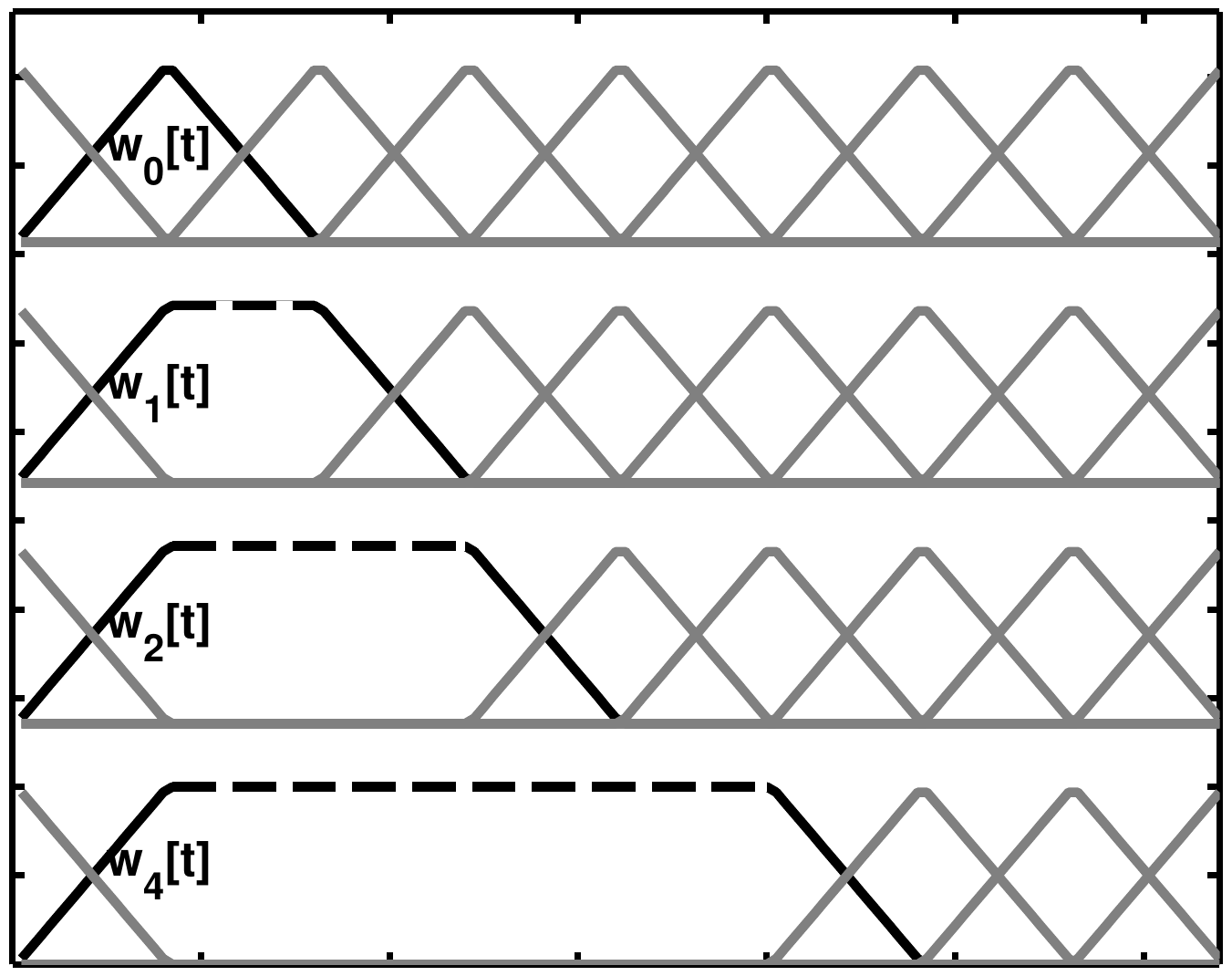}
  }\hfill
  \subfigure[\label{fig:adaptationExample}Translates of windows from a Gabor system (top) and superposition system (bottom).]{
  \includegraphics[width=.28\textwidth]{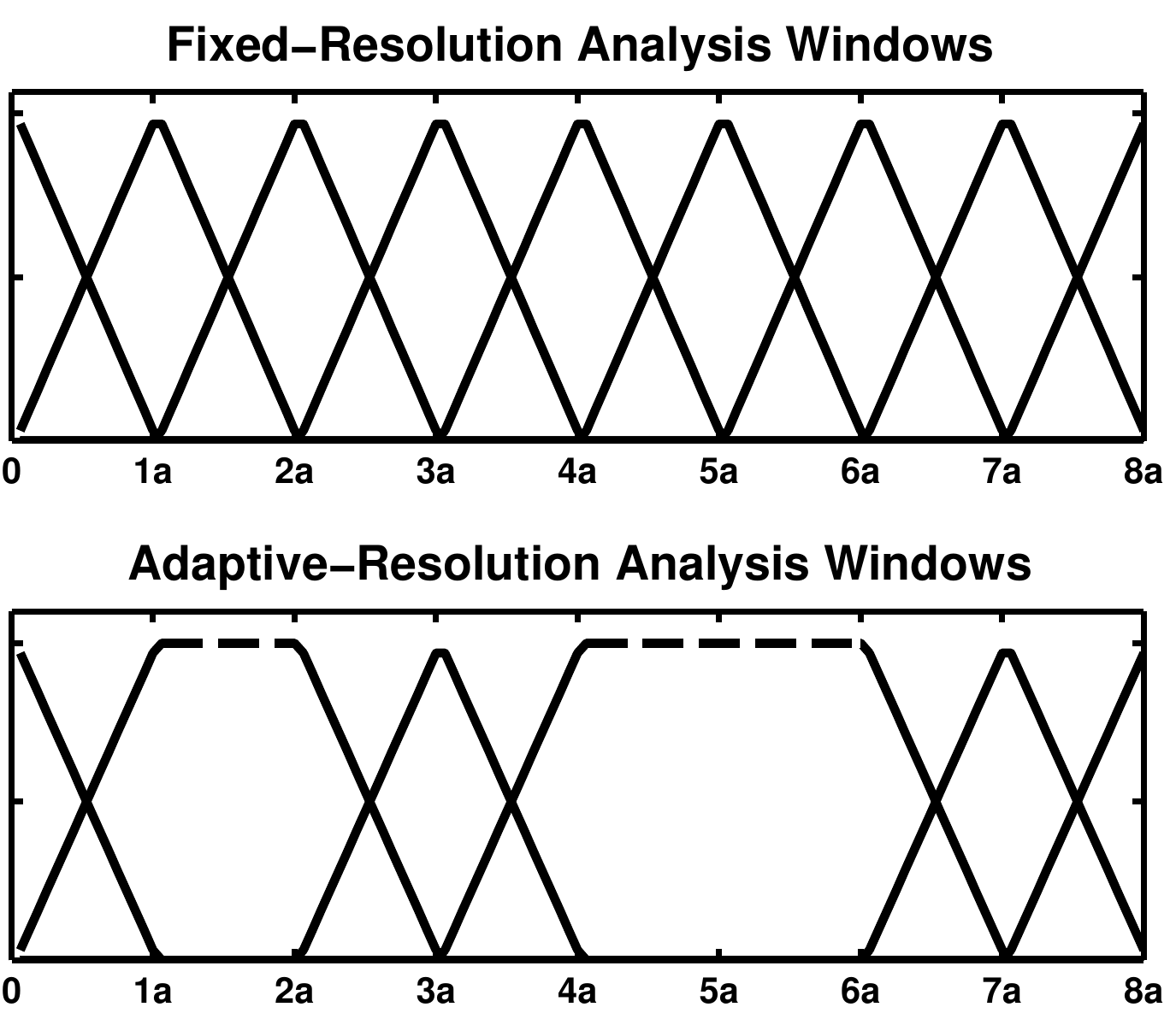}
  }\hfill
  \subfigure[\label{fig:tflattices}Corresponding time-frequency lattices associated with the systems of Fig.~\ref{fig:adaptationExample}.]{
  \includegraphics[width=.32\textwidth]{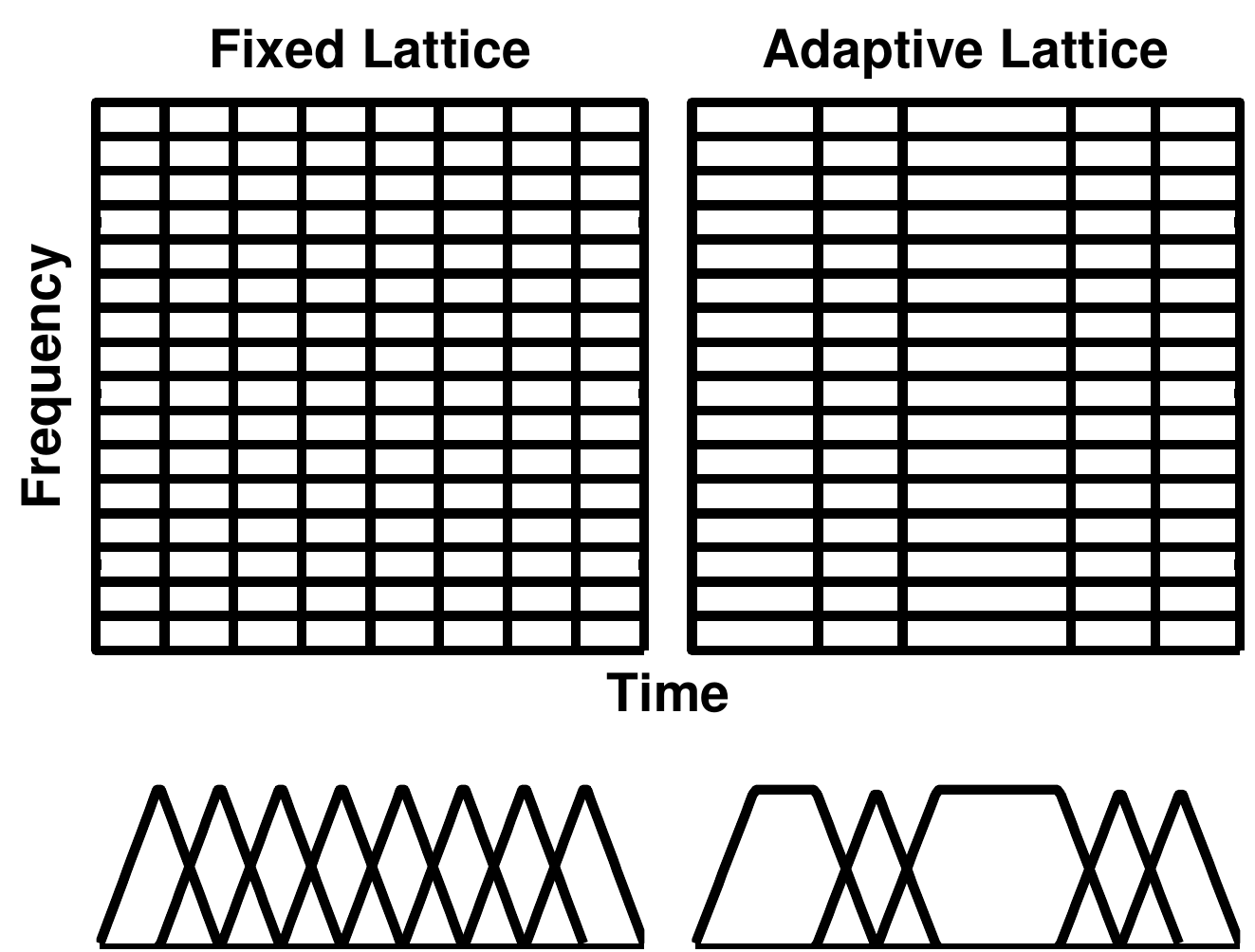}
  }\hfill
  \caption{\label{fig:adaptationExampleBoth}An example realization of a superposition system realized via two and then three window merges.}
\end{figure*}

We conclude by using the arguments of Remarks~\ref{rem:windowlength} and \ref{rem:diagFrameOpr} to establish a result required for our subsequent development.
\begin{lem}
    \label{lem:helperLemma}
    Fix any $w \in \mathbb{C}^L$ and $a,b$ such that $\mathscr{G}(w,a,b)$ is a frame for $\mathbb{C}^L$, with $M = L/b$.  Then for any integral $M' \geq \operatorname{len}(w)$, the Gabor system $\mathscr{G}(w, a, L/M')$ is also a frame for $\mathbb{C}^L$, with diagonal frame operator and maximal lower frame bound given by $(M'/M) \min_{t \in \mathbb{Z}_L}  S[t,t]$.
\end{lem}
\begin{IEEEproof}
We must show that if the frame operator $S$ of a Gabor system $\mathscr{G}(w,a,b)$ on $\mathbb{C}^L$ is full rank, then so is the frame operator $S'$ of any system $\mathscr{G}(w,a,L/M')$. For any $M' \geq \operatorname{len}(w)$, the argument of Remark~\ref{rem:diagFrameOpr} implies that $S'$ is diagonal, with eigenvalues $S'[t,t] = M' \sum_{n=0}^{N -1} |w[t - na]|^2$; the Walnut representation of~\eqref{eq:Walnut} further implies that $S'[t,t] = (M'/M) S[t,t]$, for $M = L/b$.  As $\mathscr{G}(w,a,b)$ is a frame for $\mathbb{C}^L$, it follows that $S[t,t] > 0$.  Hence $S'[t,t] > 0$ for all $t \in \mathbb{Z}_L$, and thus $S'$ is of full rank.
\end{IEEEproof}

\begin{figure}
  \centering
  \includegraphics[width=\columnwidth]{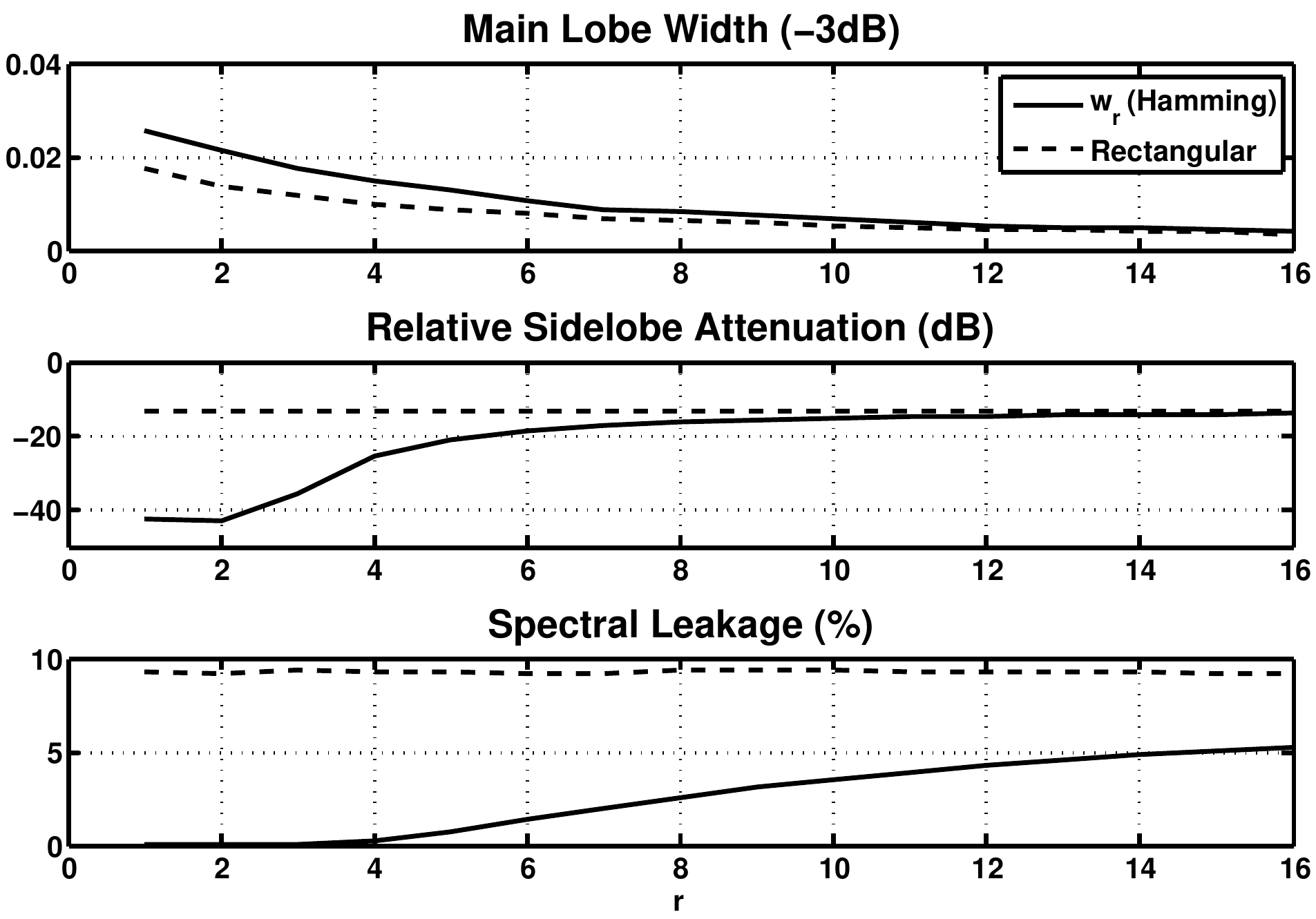}
  \caption{\label{fig:winSpectra}
  Frequency characteristics of superposition windows $w_r$ derived from Hamming windows, with $a = \operatorname{len}(w_0) / 4$ (e.g., 75\% overlap), shown relative to those of rectangular windows of length $\operatorname{len}(w_r)$.}
\end{figure}

\section{Superposition Windows}
\label{sec:superposition}

Having outlined the connections between Gabor systems and the short-time Fourier transform, we now introduce the central ingredient of our signal-adaptive time-frequency analysis framework---the \emph{superposition window} construction, illustrated in Figure~\ref{fig:windowExamples}.

\begin{defn}[Superposition Window]
\label{def:superPosition}
Fix a real, nonnegative window $w$ on $\mathbb{C}^L$ and an integer $a = L/N$, along with some $r \in \mathbb{Z}_{N}$.  We then define the  \emph{superposition window} $w_r$ to be a linear sum of $r+1$ adjacent translates of $w[t]$ as follows:
    \begin{equation}
        \label{eq:superPosWin} w_r[t] \triangleq \sum_{n=0}^{r} \mathcal{T}_{na}w[t], \quad r \in \mathbb{Z}_{N} \text{.}
    \end{equation}
\end{defn}

\begin{rem}[Fourier Transform Support]
\label{rem:supportInvariance}
Let $\widehat{w}$ denote the (discrete) Fourier transform of $w \in \mathbb{C}^L$. Linearity of~\eqref{eq:superPosWin} implies that the support of $\widehat{w_r}$ is contained within that of $\widehat{w}$, as $\operatorname{supp} (\widehat{w_r}) = \operatorname{supp} (\sum_{n=0}^r e^{-2\pi i n a (\cdot)/L} \widehat{w} ) \subseteq \operatorname{supp} (\widehat{w})$.
\end{rem}

\begin{rem}[Fourier Transform Decay]
    As $r$ increases, it is clear that $w_r$ can become more like a rectangular window (see, e.g., Fig.~\ref{fig:windowExamples}); this effect is illustrated in Fig.~\ref{fig:winSpectra} for the case of Hamming superposition windows.  Consequently, the main lobe width of $\widehat{w_r}$ shrinks, leading to improved frequency resolution relative to $\widehat{w_0}$; this main lobe resolution, however, comes at the expense of decreasing sidelobe attenuation.  Spectral leakage---a function of the window smoothness---remains superior, as does overall spectral decay for small $r$.
\end{rem}

Our subsequent construction of superposition frames employs sets of modulated superposition windows, and to this end we establish the following two energy ``conservation'' properties, proved in the appendix.

\begin{lem}[Localized Parseval Property]\label{lem:energyCons}
Fix any $w \in \mathbb{C}^L$ and an integer $M = L/b \geq \operatorname{len}(w)$.  Then
\begin{equation}
    \label{eq:energyCons}
\forall x \in \mathbb{C}^L, \,\,
\sum_{m=0}^{M-1}  \left |\langle x,\mathcal{M}_{mb}w\rangle \right |^2 = M \sum_{t=0}^{L-1} |x[t]|^2 |w[t]|^2
\text{.}
\end{equation}
\end{lem}

\begin{lem}[Superadditivity of Superposition Energy]\label{lem:emw} Let real, nonnegative superposition windows $w_p$ and $w_q$ be derived from a Gabor system $\mathscr{G}(w,a,b)$ on $\mathbb{C}^L$, and merge them to obtain a new superposition window $w_p + {w_q}' = w_p + \mathcal{T}_{(p+1)a}w_q$.  Then, \emph{if} and \emph{only if} $M = L/b \geq \operatorname{len}(w_p + {w_q}')$, the following holds for every $M_0 = L / b_0 \in \{\max(\operatorname{len}(w_p), \operatorname{len}(w_q)), \ldots, M\}$:
\begin{multline}
    \label{eq:energyLemma}
    \sum_{m=0}^{M_0-1}  \left |\langle x,\mathcal{M}_{mb_0}w_p\rangle \right |^2 + \left | \langle x,\mathcal{M}_{mb_0}{w_q}' \rangle \right |^2 \\
    \leq \sum_{m=0}^{M-1} \left|\langle x,\mathcal{M}_{mb}( w_p + {w_q}') \rangle \right|^2, \quad \forall x \in \mathbb{C}^L \text{.}
\end{multline}
\end{lem}

This superadditivity property, which is invariant to translation of $w_p + {w_q}'$, will be used in Section~\ref{sec:adaptiveFrameAnalysis} to show that adapting a Gabor frame through the superposition construction \emph{preserves} the original Gabor lower frame bound---an important consideration for numerical stability.

\section{Construction of Superposition Systems}
\label{sec:adaptFramework}

We now describe how to employ the superposition windows of Section~\ref{sec:superposition} above to create a signal-adaptive analysis framework. Let $\mathscr{G}(w, a, b)$ represent a Gabor system, which induces a short-time Fourier transform on $\mathbb{C}^L$ according to Definition~\ref{def:STFT}. Beginning with $\mathscr{G}(w, a, b)$, we then form a signal-dependent, variable-resolution STFT by using the superposition sum of~\eqref{eq:superPosWin} to adaptively merge neighboring translates from the set $\{ \mathcal{T}_{n{a}}w, n \in \mathbb{Z}_{N} \}$.  Later we will demonstrate how this signal-adaptive analysis can be coupled with a variety of different algorithms; we begin, however, by studying the general set of \emph{superposition systems} independently of any algorithmic construction.  To this end, we introduce the notion of \emph{ordered partition functions} as a means of indexing arbitrary sets of superposition windows, and then extend these to yield a full time-frequency analysis.

\subsection{Ordered Partition Functions}
\label{sec:timeLattices}

Observe that exactly $2^{N-1}$ distinct sets of variable-length superposition windows may be derived by merging window translates from a given Gabor system $\mathscr{G}(w,a,\cdot)$.  As a means of indexing these sets, the following ``stick-breaking'' analogy is helpful.  Consider a ``stick'' composed of $N$ ordered, unit-length pieces, representing elements of the set $\{ \mathcal{T}_{n{a}}w, n \in \mathbb{Z}_{N} \}$.  Merging adjacent windows in this set can be thought of as fusing neighboring pieces of the stick.  Each stick partition thus induces an ordered partition of the set $\{1,2,\ldots, N\}$, with each piece uniquely identified by an initial index and length, and we may formalize this analogy as follows.

\begin{defn}[Ordered Partition Functions]
   \label{def:stickBreaking}
   We call any $\widetilde{I}: \mathbb{Z}_{N} \times \mathbb{Z}_{N} \to \{0,1\}$ an \emph{ordered partition function} if it is not identically zero, and satisfies the following three properties:
   \begin{enumerate}
        \item \label{Iprop:one} Each piece of the stick is distinct:
        \begin{equation*}
            \widetilde{I}[n,r] = 1 \, \Rightarrow \, \widetilde{I}[n,r'] = 0 \,\,\, \forall \, r' \neq r, r' \in \mathbb{Z}_{N}  \text{.}
        \end{equation*}

        \item \label{Iprop:two} The length of each piece is denoted by $r+1$:
        \begin{equation*}
            \widetilde{I}[n,r] = 1 \, \Rightarrow \, \widetilde{I}[n',r'] = 0 \,\, \text{on} \,\, \{n+1, \ldots, n+r\} \times \mathbb{Z}_N \text{.}
        \end{equation*}

        \item \label{Iprop:three} All pieces of the stick are accounted for:
        \begin{equation*}
            \widetilde{I}[n,r] = 1 \, \Rightarrow \, \widetilde{I}[n+r,r'] = 1 \,\, \text{for exactly one} \,\, r' \in \mathbb{Z}_N \text{.}
        \end{equation*}
   \end{enumerate}
\end{defn}
Definition~\ref{def:stickBreaking} clearly implies that the ``length'' of the stick remains unchanged:
\begin{equation}
    \label{eq:Iprop_sumNR}
    \sum_{n=0}^{N-1}\sum_{r=0}^{N-1} \widetilde{I}[n,r] \,(r+1) = N \text{;}
\end{equation}
moreover, each ordered partition function $\widetilde{I}$ can be associated with a set of translated superposition windows, which includes $T_{na}w_r$ whenever $\widetilde{I}[n,r] = 1$.
The following examples of ordered partition functions are illustrated in Fig.~\ref{fig:adaptationExample}.

\begin{exmp}[$N$-Part Partition]\label{ex:uniform}
   The ordered partition function associated to the top panel of Fig.~\ref{fig:adaptationExample} is
   \begin{equation*}
       \widetilde{I}[\cdot,r] \triangleq
       \begin{cases}
       1 & \text{if $r=0$,} \\
       0 & \text{otherwise.}
       \end{cases}
   \end{equation*}
   This clearly recovers the window translates of any Gabor system $\mathscr{G}(w, a,\cdot)$. Note that in accordance with~\eqref{eq:Iprop_sumNR}, we have
   that $\sum_{n=0}^{N-1}\sum_{r=0}^{N-1} \widetilde{I}[n,r]\,(r+1) = \sum_{n=0}^{N-1}\widetilde{I}[n,0] = N = 8$.
\end{exmp}

\begin{exmp}[$(N-3)$-Part Partition]
   The ordered partition function associated to the bottom panel of Fig.~\ref{fig:adaptationExample} is
   \begin{equation*}
       \widetilde{I}[n,r] =   \begin{cases}
                       \widetilde{I}[2,0] = \widetilde{I}[6,0] = \widetilde{I}[7,0] = 1 & \text{not merged,}\\
                       \widetilde{I}[0,1] = \widetilde{I}[3,2] = 1 & \text{merged.}
                   \end{cases}
   \end{equation*}
   Note again that in accordance with~\eqref{eq:Iprop_sumNR}, we have that $\sum_{n=0}^{N-1}\sum_{r=0}^{N-1}\widetilde{I}[n,r] (r+1) = 3 + 2 + 3 = 8$.
\end{exmp}

\subsection{Superposition Systems}
\label{sec:superPosSystems}

We now employ the above construction to arrive at a variable-resolution time-frequency analysis via superposition windows.  To this end, let the set $\mathscr{F}$ be defined as a function of any Gabor system $\mathscr{G}(w,a,b)$ on $\mathbb{C}^L$ as follows:
\begin{equation*}
    \mathscr{F} \triangleq \bigcup_{r \in \mathbb{Z}_N}\mathscr{G}(w_r,a,b_L) \text{,}
\end{equation*}
where the frequency lattice $b_L \mathbb{Z}$ encompasses all possible Gabor systems on $\mathbb{C}^L$ for a fixed choice of integral $M$:
\begin{equation*}\label{eq:bLDef}
    b_L \mathbb{Z}; \,
    M_L \triangleq \operatorname{lcm}(\{1,2 \ldots, \max(L,M)\}), \,
    b_L \triangleq L/M_L \text{.}\!
\end{equation*}

Elements of $\mathscr{F}$ may then be defined in analogy to~\eqref{eq:GaborFrame} as
\begin{equation*}
 \phi_{m,n,r} \triangleq \mathcal{M}_{mb_L} \mathcal{T}_{na}w_r = \textstyle\mathcal{M}_{mb_L} \!\!\left ( \sum_{n'=0}^{r} \mathcal{T}_{(n'+n)a}w \right ) \!\text{,}
\end{equation*}
and in turn give rise to \emph{superposition systems}, defined as appropriately chosen subsets of $\mathscr{F}$.
\begin{defn}[Superposition Systems and Admissibility]\label{def:superPosSys}
   Fix an ordered partition function $\widetilde{I}[n,r]$ and a function $M[n,r]: \mathbb{Z}_N \times \mathbb{Z}_N \rightarrow \mathbb{Z}_{M_L}$.
   We call any $I[m,n,r]: \mathbb{Z}_{M_L} \times \mathbb{Z}_{N} \times \mathbb{Z}_{N} \to \{0,1\}$ an \emph{admissible selection function} on $\mathscr{F} = \cup_r\mathscr{G}(w_r,a,b_L) $ if it satisfies the following two properties:
   \begin{align}
        \nonumber
             I[0, n, r] & \!=\! \widetilde{I}[n,r] \quad \forall\, n,r \in \mathbb{Z}_N \times \mathbb{Z}_N \text{,} \\
        \label{Iprop:frequency}
             I[0, n, r] & \!=\! 1 \, \Rightarrow I[{\textstyle\frac{M_L}{M[n,r]}}m,n,r] = 1, \,\,  \forall\, m \in Z_{M[n,r]} \text{.}
   \end{align}
   Furthermore, we call the induced set of elements a \emph{superposition system} $\mathscr{F}(I)$:
    \begin{equation*}
        \label{eq:F_I}
        \phi_{m,n,r} \in \mathscr{F}(I) \Leftrightarrow I[m,n,r] = 1 \text{.}
   \end{equation*}
\end{defn}

It follows from~\eqref{Iprop:frequency} that the first $M[n,r]$ modulates of each selected superposition window are included in $\mathscr{F}(I)$, and thus we later suppress the dependence of $I$ on frequency bin index $m$ when possible, by abbreviating $I[\cdot,n,r]$ as $I[n,r]$.

\section{Superposition Frames: Main Results}
\label{sec:adaptiveFrameAnalysis}

Starting from a Gabor system $\mathscr{G}(w, a, b)$, we see that any superposition system $\mathscr{F}(I) \subset \cup_r\mathscr{G}(w_r,a,b_L)$ effectively yields a ``variable-resolution'' subsampled short-time Fourier transform, defined for all $m \in \mathbb{Z}_{M_L}$ and $n,r \in \mathbb{Z}_{N}$ as
\begin{equation}
   \label{eq:adaptInnerProduct}
   X[m,n,r] \triangleq
   \begin{cases}
         \langle x, \phi_{m,n,r} \rangle & \text{if $\phi_{m,n,r} \in \mathscr{F}(I)$,} \\
         0 & \text{otherwise.}
   \end{cases}
\end{equation}
Consequently, we now establish conditions under which superposition systems $\mathscr{F}(I)$ form frames for $\mathbb{C}^L$, in analogy to the relation between a Gabor frame and its corresponding fixed-resolution short-time Fourier transform.

\subsection{General Case: Sufficiency}

To begin our analysis, consider the case of an admissible selection function $I[m,n,r]$ for which  $M[n,r] = M_{\textrm{g}}$ for all $n,r \in \mathbb{Z}_N$, corresponding to the notion of a \emph{global} frequency lattice of arbitrary resolution: $b_{\textrm{g}} \mathbb{Z}$ with $b_{\textrm{g}} = L / M_{\textrm{g}}$.  Our first result, proved in the appendix, ensures that the induced superposition system $\mathscr{F}(I)$ is a frame for $\mathbb{C}^L$ if the following test condition holds.

\begin{thm}[Sufficiency$\!$ Condition,$\!$ Superposition$\!$ Frames]\label{thm:testCond}
Fix a Gabor system $\mathscr{G}(w, a, \cdot)$ on $\mathbb{C}^L$, with $N\!=\!L/a$, $w$ real and nonnegative, and define for $s,t \!\in\! \mathbb{Z}_L$, $n,r \!\in\! \mathbb{Z}_{N}$, the term
\begin{equation*}
  \beta_{nr}(s,t) \triangleq w_r[t-na] w_r[t-na-s] \text{.}
\end{equation*}
Let $I[n,r]$ be \emph{any} admissible selection function for which $M[n,r] = M_{\textrm{g}}$ for some $M_{\textrm{g}} \in \{1,2,\ldots, L\}$. Then, for index term $k \in \{ \lceil (t-(L-1))/M_{\textrm{g}}  \rceil, \ldots, \lfloor t/M_{\textrm{g}} \rfloor \}$, the condition
\begin{equation}
  \label{eq:testCondition}
  \forall t \in \mathbb{Z}_L, \,
  \sum_{n=0}^{N-1} \sum_{r=0}^{N-1} I[n,r] \Big( \beta_{nr}(0,t) - \sum_{k \neq 0} \beta_{nr}(kM_{\textrm{g}},t) \Big) > 0
\end{equation}
is sufficient to guarantee that the superposition system $\mathscr{F}(I) \subset \cup_{r} \mathscr{G}(w_r,a,b_L)$ is a frame for $\mathbb{C}^L$.
\end{thm}

Satisfying the criterion of~\eqref{eq:testCondition} implies that the underlying frame operator is \emph{strictly diagonally dominant}---a sufficient condition for strict positive definiteness. This is a popular criterion in the literature (see, e.g.,~\cite[Corollary~6]{Qiu1995},~\cite[Theorem~8.4.4]{christensen2003ifa}) and, as can be seen from~\eqref{eq:testCondition}, takes a particularly simple form in the superposition setting.

\subsection{Superposition Frames and Frame Bounds}

In Theorem~\ref{thm:testCond} above, we considered a general class of superposition frames associated with an arbitrary frequency lattice $b_{\textrm{g}} \mathbb{Z}$. In Theorems~\ref{thm:adaptiveFrames} and~\ref{thm:superPosFrameBounds} below, we study two distinct classes of superposition systems using non-uniform (local) and uniform (global) modulation structures defined as follows.
\begin{defn}[Admissible Selection Functions ${I^l}$ and ${I^{\textrm{g}}}$]
    \label{defn:localGlobalI}
    Fix a Gabor system $\mathscr{G}(w,a,L/M)$, associate to it any ordered partition function $\widetilde{I}[n,r]$, and define
    \begin{align}
        \label{eq:MLocal}
            M_r & \triangleq \,\,  \max \, \left( \operatorname{len}(w_r),M \right); \quad b_r \triangleq L/M_r \text{,} \\
            \label{eq:MGlobal}
            M_{\textrm{g}} & \triangleq \!\! \max_{r \,:\, \widetilde{I}[\cdot, r] \,=\, 1 } M_r; \qquad \qquad b_{\textrm{g}} \triangleq L/M_{\textrm{g}}  \text{.}
    \end{align}
    These quantities induce, via $M[n,r] = M[\cdot,r] = M_r$ or $M[n,r] = M_{\textrm{g}}$ constant, respective classes of \emph{local} and \emph{global} admissible selection functions $I^l[m,n,r]$ and $I^{\textrm{g}}[m,n,r]$.  Note that when the admissible selection functions $I_0$ or $I_{N-1}$ occur, the constants of~\eqref{eq:MLocal} and~\eqref{eq:MGlobal} are equal.
\end{defn}

We now show that superposition systems $\mathscr{F}(I^{\textrm{g}})$ and $\mathscr{F}(I^l)$ are frames for $\mathbb{C}^L$. Later, we will verify that such frames admit \emph{diagonal} frame operators. This special structure leads not only to fast reconstruction algorithms, but also to the preservation of lower frame bounds.

\begin{thm}[Local and Global Superposition Frames]
\label{thm:adaptiveFrames}
    Let $\mathscr{G}(w, a, b)$ be a Gabor frame for $\mathbb{C}^L$, with $w$ real and nonnegative. Then for any choice of admissible selection functions $I^l$ and $I^{\textrm{g}}$, the local and global superposition systems $\mathscr{F}(I^l)$ and $\mathscr{F}(I^{\textrm{g}})$ are also frames for $\mathbb{C}^L$.
\end{thm}

\begin{IEEEproof}
As our finite-dimensional setting implies the existence of an upper frame bound for any admissible $I$, only the existence of a lower frame bound need be established. The proof proceeds via Lemma~\ref{lem:emw} and an iterative argument.

To begin, consider the admissible selection function $I^{\textrm{g}}_0[m,n,r] = I^l_0[m,n,r]$ induced by an $N$-part ordered partition, which is associated to the event that \textit{no merging} of windows in $\mathscr{G}(w,a,b)$ occurs, and hence $\mathscr{F}(I^{\textrm{g}}_0) = \mathscr{G}(w,a,b_{\textrm{g}})$, with $b_{\textrm{g}} = L/M_{\textrm{g}}$ and $M_{\textrm{g}} = \max(\operatorname{len}(w),M)$ according to~\eqref{eq:MGlobal}. Lemma~\ref{lem:helperLemma} then ensures that $\mathscr{F}(I^{\textrm{g}}_0) = \mathscr{F}(I^l_0)$ is a frame for $\mathbb{C}^L$, with maximal lower frame bound $M_\textrm{g} = \min_{t\in\mathbb{Z}} \sum_{n=0}^{N-1} |\mathcal{T}_{na}w[t]|^2 > 0$.

Next consider any admissible selection function $I^{\textrm{g}}_1[m,n,r]$ induced by an $(N-1)$-part ordered partition, corresponding to the case that \emph{exactly one} pair of windows $w \equiv w_0$ from the initial Gabor frame $\mathscr{G}(w, a, b)$ is merged via the superposition sum of~\eqref{eq:superPosWin}. In this case, there exists one $n^* \in \mathbb{Z}_{N}$ such that $I^{\textrm{g}}_1[0,n^*,1] = 1$, and so~\eqref{Iprop:frequency} implies that $\mathscr{F}(I^{\textrm{g}}_1)$ contains the elements $\{\mathcal{M}_{mb_{\textrm{g}}}\mathcal{T}_{n^*a}w_{1}: m \in \mathbb{Z}_{M_{\textrm{g}}}\}$.  Each of these elements can in turn be decomposed into the following sum:
\begin{equation}
    \label{eq:iterStep1}
  \mathcal{M}_{mb_{\textrm{g}}}\mathcal{T}_{n^*a}w_1 = \mathcal{M}_{mb_{\textrm{g}}}\mathcal{T}_{n^*a}w_0 + \mathcal{M}_{mb_{\textrm{g}}}\mathcal{T}_{(n^*+1)a}w_0
  \text{.}
\end{equation}

Since~\eqref{eq:MGlobal} implies $M_{\textrm{g}} \geq \operatorname{len}(w_1) = \operatorname{len}(w_0 + \mathcal{T}_{a}w_0)$, we obtain by~\eqref{eq:iterStep1} and the superadditivity property of Lemma~\ref{lem:emw}:
\begin{multline}\label{eq:sumIneq}
    \sum_{m=0}^{M_{\textrm{g}}-1}  | \langle x, \mathcal{M}_{mb_{\textrm{g}}}\mathcal{T}_{n^*a}w_{1} \rangle  |^2
    \geq  \sum_{m=0}^{M_{\textrm{g}}-1} |\langle x, \mathcal{M}_{mb_{\textrm{g}}}\mathcal{T}_{n^*a}w_0 \rangle |^2
    \\
    +  \sum_{m=0}^{M_{\textrm{g}}-1} |\langle x, \mathcal{M}_{mb_{\textrm{g}}}\mathcal{T}_{(n^*+1)a}w_0 \rangle |^2
    , \quad \forall x \in \mathbb{C}^L \text{.}
\end{multline}

Next, noting that by the decomposition of~\eqref{eq:iterStep1}, we have
\begin{multline}\label{eq:breakPoint}
 \mathscr{F}(I^{\textrm{g}}_1) = \mathscr{G}(w,a,b_{\textrm{g}}) \cup \{\mathcal{M}_{mb_{\textrm{g}}}\mathcal{T}_{n^*a}w_1\}
 \\
 \setminus \left ( \{\mathcal{M}_{mb_{\textrm{g}}}\mathcal{T}_{n^*a}w_{0}\} \cup \{\mathcal{M}_{mb_{\textrm{g}}}\mathcal{T}_{(n^*+1)a}w_{0}\} \right )
 \text{,}
\end{multline}
we see that~\eqref{eq:sumIneq} and~\eqref{eq:breakPoint} together imply that for all $x \in \mathbb{C}^L$,
\begin{equation}
    \label{eq:iter1Inequality}
    \sum_{\phi  \in \mathscr{F}(I^{\textrm{g}}_1)} |\langle x,\phi_{m,n,r} \rangle |^2
    \geq \!\!\! \sum_{\phi \in \mathscr{G}(w,a,b_{\textrm{g}})}  |\langle x,\phi_{m,n} \rangle|^2 \text{.}
\end{equation}
Since $\mathscr{G}(w,a,b_{\textrm{g}}) = \mathscr{F}(I^{\textrm{g}}_0)$ is a frame, the existence of a lower frame bound for $\mathscr{F}(I^{\textrm{g}}_1)$ is guaranteed by~\eqref{eq:iter1Inequality}.

Now consider the general case in which $\mathscr{F}(I^{\textrm{g}})$ contains multiple merges. Since our construction ensures that all admissible selection functions can be obtained by iterative partitioning in the manner above,  we can always recover the inequality of~\eqref{eq:iter1Inequality} for $\mathscr{F}(I^{\textrm{g}})$ and any  $b_{\textrm{g}}$ according to~\eqref{eq:MGlobal}, by linearity of superposition and repeated application of Lemma~\ref{lem:emw}.

A similar iterative argument holds for $\mathscr{F}(I^l_1)$, with $b_{\textrm{g}}$ replaced by $b_0$ from~\eqref{eq:MLocal}. In place of~\eqref{eq:breakPoint} we obtain
    \begin{multline}\label{eq:breakPointLocal}
            \mathscr{F}(I^l_1) = \left (\mathscr{G}(w,a,b_0) \cup \{\mathcal{M}_{mb_1}\mathcal{T}_{n^*a}w_1\} \right )
            \\ \setminus \left ( \{\mathcal{M}_{mb_0}\mathcal{T}_{n^*a}w_{0}\} \cup \{\mathcal{M}_{mb_0}\mathcal{T}_{(n^*+1)a}w_{0}\} \right ) \text{,}
    \end{multline}
    with $b_1=L/M_1$ according to~\eqref{eq:MLocal}.  Note that
    \begin{align*}
    M_1 & \geq \operatorname{len}(w_1) = \operatorname{len}(w_0 + \mathcal{T}_{a}w_0)\text{,}
    \\ M_0 & \geq \operatorname{len}(w_0) \text{,~~and } M_1 \geq M_0
    \text{;}
    \end{align*}
    thus, we may apply Lemma~\ref{lem:emw} to the latter three terms of~\eqref{eq:breakPointLocal}, yielding the required result for $\mathscr{F}(I^l_1)$: for all
    $x \in \mathbb{C}^L$,
    \begin{equation*}
        \sum_{\phi  \in \mathscr{F}(I^l_1)} |\langle x,\phi_{m,n,r} \rangle |^2
        \geq \!\!\! \sum_{\phi \in \mathscr{G}(w,a,b_0)}  |\langle x,\phi_{m,n,r} \rangle|^2 \text{.}
    \end{equation*}
    As above, the proof for general $\mathscr{F}(I^l)$ then follows.
\end{IEEEproof}

Thus we see that for every Gabor system on $\mathbb{C}^L$ and any associated ordered partition function, setting $M[n,r]$ in accordance with $M_r$ or $M_{\textrm{g}}$ will yield local or global superposition frames. Moreover, as we detail later, the iterative arguments employed above suggest precise algorithmic constructions.

We now proceed to establish the important property that, for all local and global ${I[m,n,r]}$ of Definition~\ref{defn:localGlobalI}, \emph{superposition frames preserve lower frame bounds}, thus ensuring numerical stability of the resultant representation.  The following result also formulates the corresponding minimax-optimal superposition frame bounds.

\begin{thm}[Superposition Frame Bound Properties]
  \label{thm:superPosFrameBounds}
  Let $\mathscr{G}(w, a, b)$ be a frame for $\mathbb{C}^L$, with associated maximal lower frame bound $A > 0$.  Then for any admissible $I^{l}$ and $I^{\textrm{g}}$:
  \begin{enumerate}

  \item\label{claim:frameBoundPres} The quantity $A$ remains a valid lower frame bound for both $\mathscr{F}(I^l)$ and $\mathscr{F}(I^{\textrm{g}})$.

  \item The \emph{minimum} maximal lower superposition frame bound
over all admissible $I^{l}$ and $I^{\textrm{g}}$ is
      \begin{equation*}
      A_{\mathrm{opt}} = \frac{L}{b^{\mathrm{max}}} \cdot \min_{t \in
\mathbb{Z}_L} \left( \sum_{n=0}^{N -1} |\mathcal{T}_{na}w[t]|^2
\right) \text{,}
      \end{equation*}
      with $b^{\text{max}} =\min(b,L/\operatorname{len}(w))$. It is
attained in the absence of merging: $\mathscr{F}(I^l) =
\mathscr{F}(I^{\textrm{g}})  = \mathscr{G}(w,a,b^{\mathrm{max}})$.

  \item The \emph{maximum} minimal upper superposition frame bound
over all admissible $I^{l}$ and $I^{\textrm{g}}$ is
      \begin{equation*}
       B_{\mathrm{opt}} =
         \frac{L}{b^{\mathrm{min}}} \cdot \max_{t \in \mathbb{Z}_L}
\left( \Big|\sum_{n=0}^{N -1} \mathcal{T}_{na}w[t] \Big|^2 \right)
\text{,}
      \end{equation*}
      with $b^{\mathrm{min}} = \min(b,1)$.  It is attained when all
translates of $w$ have been merged: $\mathscr{F}(I^l) =
\mathscr{F}(I^{\textrm{g}}) =
\mathscr{G}(w_{N-1},L,b^{\mathrm{min}})$, with  $w_{N-1} =
\sum_{n=0}^{N -1} \mathcal{T}_{na}w[t]$.
  \end{enumerate}
\end{thm}

\begin{IEEEproof}
Let $S$ be the frame operator associated to $\mathscr{G}(w, a, b)$, with smallest eigenvalue $A$, and observe that $\min_{t \in \mathbb{Z}_L} S[t,t] \geq A$ by the Schur-Horn convexity theorem. Now, for any admissible $I^{\textrm{g}}$ or $I^l$, the proof of Theorem~\ref{thm:adaptiveFrames} shows that the maximal lower frame bound of $\mathscr{F}(I^{\textrm{g}})$ or $\mathscr{F}(I^{l})$ is bounded from below by that of some $\mathscr{G}(w, a, b')$, with $b'$ denoting respectively $b_{\textrm{g}}$ or $b_0$. Therefore,~\ref{claim:frameBoundPres}) will follow if we can show the maximal lower frame bound of $\mathscr{G}(w, a, b')$ to be no less than $\min_{t \in \mathbb{Z}_L} S[t,t]$.  To do so, note that Lemma~\ref{lem:helperLemma} implies that the frame operator $S'$ of $\mathscr{G}(w, a, b')$ is diagonal, with smallest eigenvalue $\min_{t \in \mathbb{Z}_L} S'[t,t] = (M'/M) \min_{t \in \mathbb{Z}_L} S[t,t]$;~\eqref{eq:MLocal} and~\eqref{eq:MGlobal} then yield $M' \geq M$.

Next recall that superposition frames are induced from a Gabor frame $\mathscr{G}(w, a, b)$ by merging $n \in \mathbb{Z}_N$ neighboring window translates $\mathcal{T}_{na}w$.  By the above argument, $\min_{t \in \mathbb{Z}_L} S'[t,t]$ is itself a maximal lower frame bound for the case $\mathscr{G}(w,a,b')$ attained whenever \emph{no} window translates are merged; likewise, the merging of \emph{all} translates yields $\mathscr{G}(w_{N-1},L,b'')$ for some unique $b''$, with minimal upper frame bound $\max_{t \in \mathbb{Z}_L} S''[t,t]$.

To show that these cases are in fact extremal as claimed, we appeal to the same iterative argument used to prove Theorem~\ref{thm:adaptiveFrames}.  There, the superadditivity property of Lemma~\ref{lem:emw} was invoked to show that for any admissible $I^{\textrm{g}}$ or $I^l$, merging superposition windows \emph{cannot decrease} the overall energy of the resultant frame coefficients. Thus, the case of $\mathscr{G}(w,a,b')$ considered above represents attainment of the \emph{minimum} maximal lower superposition frame bound, and moreover $M' = \max(M,\operatorname{len}(w))$ via~\eqref{eq:MLocal} and~\eqref{eq:MGlobal}.  Likewise, $\mathscr{G}(w_{N-1},L,b'')$ yields the \emph{maximum} minimal upper superposition frame bound, with $M'' = \max(M,\operatorname{len}(w_{N-1}))$. Lemma~\ref{lem:energyCons} then establishes the bound directly:
\begin{align*}
 &  
 \sum_{\phi \in \mathscr{G}(w_{N-1},L,b'')} \!\!\!\! |\langle x,
\phi_{m,n,r} \rangle|^2 = \textstyle \sum_{m=0}^{M''-1} | \langle
x,\mathcal{M}_{mb''}w_{N-1} \rangle|^2
\\ &  = M'' \textstyle \sum_{t=0}^{L-1} |x[t]|^2 |w_{N-1}[t]|^2
 \leq  \|x\|^2 M'' \displaystyle \max_{t \in \mathbb{Z}_L} |w_{N-1}[t]|^2
\text{,}
\end{align*}
and the proof is completed by noting that as $\mathscr{G}(w,a,b)$ is assumed a Gabor frame for $\mathbb{C}^L$, the covering condition of Remark~\ref{rem:covCond} implies that $|\operatorname{supp}(w_{N-1})| = \operatorname{len}(w_{N-1}) = L$, and hence $M'' = \max(M,L)$ as claimed.
\end{IEEEproof}

\section{Fast Reconstruction via Superposition Frames}
\label{sec:synthesis}

We now show how the special structure of our superposition construction gives rise to a number of efficient reconstruction procedures.  For any signal of interest $x \in \mathbb{C}^L$, the superposition frame analysis coefficients  $X[m,n,r] = I[m,n,r]  \langle x, \phi_{m,n,r} \rangle $ can be computed via fast Fourier transform (FFT) once an admissible selection function has been specified.  Superposition frames also enable \emph{fast} (FFT-based) \emph{reconstruction} from the corresponding analysis coefficients, in contrast to the general case of $\mathcal{O}(L^3)$ complexity for frame-based reconstruction via inversion of the frame operator.

We first provide a fast constant-overlap-add reconstruction method, which obviates the need for canonical dual frames. We next show that reconstruction via the canonical dual can also proceed by way of a pointwise modification of each superposition window $\phi_{0,n,r}$, followed by the application of FFTs, as in the case of general nonstationary Gabor frames~\cite{JailletDorfler09}.  Third, we show that in settings reminiscent of lapped orthogonal transforms, calculation of canonical dual windows is possible \emph{independently} of any $I^\textrm{g}$---in contrast to the typical signal-adaptive setting, where the structure of the frame operator is a function of the instantiated signal adaptation.  Last, we compare the computational complexity of these procedures.

\subsection{Reconstruction via the Constant Overlap-Add Method}
\label{sec:Ola}

The classical ``overlap-add'' approach to signal reconstruction from short-time Fourier coefficients proceeds as follows~\cite{QuatieriChapter}. Recall the covering condition of~\eqref{eq:covCond}
which is necessary for a Gabor system $\mathscr{G}(w,a,b)$ to form a frame for $\mathbb{C}^L$, and also sufficient if $M = L/b \geq \operatorname{len}(w)$. Clearly, this covering condition holds if translates $\{ \mathcal{T}_{na}w : n \in \mathbb{Z}_N\}$ form a partition of unity on $\mathbb{C}^L$ (see, e.g., the top panel of Fig.~\ref{fig:adaptationExample}), and to this end we obtain the following definition, long popular in the signal processing literature.
\begin{defn}[Constant Overlap-Add Window Constraint]\label{defn:olaConstraint}
    Fix a Gabor system $\mathscr{G}(w, a, \cdot)$ on $\mathbb{C}^L$.  Then, noting the discrete Fourier transform evaluation $\widehat{w}[0] = \sum_{t=0}^{L-1}w[t]$, the window $w$ is said to satisfy the \emph{constant overlap-add} constraint if
    \begin{equation}
        \label{eq:olaConstraint}
        \forall \, t \in \mathbb{Z}_L, \,\,\, \sum_{n = 0}^{N-1}  w[t-na] = \frac{\widehat{w}[0]}{a} \text{.}
    \end{equation}
\end{defn}

To clarify the role of this overlap-add constraint in fast reconstruction, consider a Gabor frame $\mathscr{G}(w,a,b)$ on $\mathbb{C}^L$ for which $w$ satisfies~\eqref{eq:olaConstraint}, with $M = L/b$ chosen such that $M \geq \operatorname{len}(w)$.  The associated short-time analysis coefficients $\{ X[m,n] \, : \, m \in \mathbb{Z}_{M}, n \in \mathbb{Z}_N \}$ are obtained as inner products of any $x \in \mathbb{C}^L$ according to~\eqref{eq:STFT}, and it is easy to show that
\begin{equation}
    \label{eq:olaReconstruction}
    x[t] = \frac{a}{\widehat{w}[0]}\sum_{n=0}^{N-1} \left( \frac{1}{M} \sum_{m=0}^{M-1} X[m,n]e^{2\pi imbt/L} \right) \text{.}
\end{equation}
The constraint thus admits a reconstruction procedure based on the overlapping additions of a sequence of discrete Fourier transforms on $\mathbb{C}^M$.

\begin{rem}[Superposition$\!$ Windows$\!$ Preserve$\!$ Overlap-Add]\label{rem:olaPres}
By their \emph{linear} construction, superposition windows $w_r$ preserve the constant overlap-add constraint of Definition~\ref{defn:olaConstraint} for any Gabor system $\mathscr{G}(w_r,a,\cdot)$. To see this, note that the discrete Poisson summation formula on $\mathbb{C}^L$, for $N = L/a$, is given by $
    \sum_{n=0}^{N-1} \mathbb{I}_{0}[t-na] = a^{-1}\sum_{k = 0}^{a-1} e^{2\pi i t k N/L}$.
Applying this expression to $\{ \mathcal{T}_{na}w : n \in \mathbb{Z}_N\}$ yields the relation
\begin{equation*}
    \sum_{n = 0 }^{N-1} w[t-na] = \frac{1}{a} \sum_{k = 0}^{a-1} e^{2\pi i t k N/L} \widehat{w}[kN] \text{,}
\end{equation*}
and it follows that the constraint of~\eqref{eq:olaConstraint} holds (for a given time lattice constant $a$) if the Fourier transform $\widehat{w}$ satisfies
\begin{equation*}
         \widehat{w}[kN] = 0, \,\, \forall\, k \in \{1, \ldots, a-1\} \text{.}
    \end{equation*}
Since $\operatorname{supp} (\widehat{w_r}) \subseteq \operatorname{supp} (\widehat{w})$, in accordance with the argument of Remark~\ref{rem:supportInvariance}, it follows that if~\eqref{eq:olaConstraint} holds for a given  $\mathscr{G}(w,a,\cdot)$, then it will also hold for any $\mathscr{G}(w_r,a,\cdot)$, for $r \in \mathbb{Z}_N$.
\end{rem}

The popularity of the overlap-add constraint of Definition~\ref{defn:olaConstraint} is due in large part to its simplicity, coupled with the efficiency of evaluating~\eqref{eq:olaReconstruction}. An important property of our superposition construction is that it preserves this constraint not only for Gabor frames $\mathscr{G}(w_r,a,b)$, but also for all induced superposition frames $\mathscr{F}(I^{\textrm{g}})$ and $\mathscr{F}(I^l)$.

\begin{thm}[Superposition Frames Preserve Overlap-Add]\label{thm:GOLA}
Consider a Gabor frame $\mathscr{G}(w,a,b)$ on $\mathbb{C}^L$ satisfying the constant overlap-add constraint of~\eqref{eq:olaConstraint}. The following statements hold for any $I^{\textrm{g}}$, and also for any $I^l$, with $M_{\textrm{g}}, b_{\textrm{g}}$ replaced by $M_r, b_r$.

\begin{enumerate}
    \item\label{item:GOLA} The superposition frame $\mathscr{F}(I^{\textrm{g}})$ satisfies the following \emph{generalized} overlap-add constraint:
    \begin{equation}
        \label{eq:adaptiveOlaConstraint}
        \sum_{n=0}^{N-1} \sum_{r=0}^{N-1} I^{\textrm{g}}[n,r]\mathcal{T}_{na}w_r[t] = \frac{\widehat{w}[0]}{a} \text{.}
    \end{equation}

    \item\label{item:GOLARecons} Each $\mathscr{F}(I^{\textrm{g}})$ satisfies the overlap-add reconstruction property that, for any $x \in \mathbb{C}^L$ with corresponding frame coefficients $\{ X[m,n,r] \}$ defined by~\eqref{eq:adaptInnerProduct},
    \begin{equation}
        \label{eq:golaReconstruction}
        x[t] \!=\! \frac{a}{\widehat{w}[0]} \sum_{n=0}^{N-1} \! \sum_{r=0}^{N-1} \!\left(\! \frac{1}{M_{\textrm{g}}} \!\!\sum_{m=0}^{M_{\textrm{g}}-1} \!\!\!X[\textstyle{\frac{M_L}{M_{\textrm{g}}}}m,n,r] e^{2\pi imb_{\textrm{g}}t/L} \!\right)\!\!
        \text{.}
    \end{equation}
\end{enumerate}
\end{thm}

\begin{IEEEproof}
    First consider $I^{\textrm{g}}_0[n,r]$, the global selection function associated to the event that no windows are merged. In this case, the generalized overlap-add constraint of~\eqref{eq:adaptiveOlaConstraint} is satisfied by hypothesis, since it follows directly from~\eqref{eq:olaConstraint} that
    \begin{equation*}
        \sum_{n=0}^{N-1} \sum_{r=0}^{N-1} I^{\textrm{g}}_0[n,r]\mathcal{T}_{na}w_r[t] = \sum_{n=0}^{N-1} \mathcal{T}_{na}w[t] = \frac{\widehat{w}[0]}{a} \text{.}
    \end{equation*}
    Now consider the global selection function $I^{\textrm{g}}_1[n,r]$, corresponding to the case that \textit{exactly one} pair of windows $w \equiv w_0$ is merged.  In this case, there exists some $n^* \in \mathbb{Z}_N$ such that $I^{\textrm{g}}_1[n^*,1] = 1$. Thus, the frame $\mathscr{F}(I^{\textrm{g}}_1)$ contains the element $\{\mathcal{T}_{n^*a}w_{1}\}$, which can be decomposed according to~\eqref{eq:superPosWin} as
    \begin{equation}\label{eq:OLAsplit}
        \mathcal{T}_{n^*a}w_1 = \mathcal{T}_{n^*a}w_0 + \mathcal{T}_{(n^*+1)a}w_0 \text{.}
    \end{equation}
    Admissibility of $I^{\textrm{g}}_1$ implies that $I^{\textrm{g}}_1[n^*,1] = 1$, but that $I^{\textrm{g}}_1[n^*,0] = I^{\textrm{g}}_1[n^*\!+\!1,0] = 0$. Therefore, by~\eqref{eq:OLAsplit} we have
    \begin{align*}
        \sum_{n=0}^{N-1} \sum_{r=0}^{N-1} I^{\textrm{g}}_1[n,r]\mathcal{T}_{na}w_r[t] & = \mathcal{T}_{n^*a}w_1[t] + \!\!\!\! \!\!\!\! \!\!\! \!\!\! \!\!\! \!\!\!
        \sum_{\qquad n \in \mathbb{Z}_N \setminus \{n^*, n^*\!+\!1 \}}
        \!\!\!\! \!\!\! \!\!\! \!\!\!\! \!\!\!\! \! I^{\textrm{g}}_1[n,0] \mathcal{T}_{na}w_0[t] \\
        & = \! \sum_{n=0}^{N-1} \mathcal{T}_{na}w_0[t] = \frac{\widehat{w}[0]}{a} \text{,}
    \end{align*}
    and we see that~\eqref{eq:adaptiveOlaConstraint} holds for $\mathscr{F}(I^{\textrm{g}}_1)$. Naturally, the selection function $I^{\textrm{g}}$ may index many merged windows---not just one, as in the case of $I^{\textrm{g}}_1$. However, repeated application of the above argument shows that $\mathscr{F}(I^{\textrm{g}})$ satisfies~\eqref{eq:adaptiveOlaConstraint} for any $I^{\textrm{g}}$.

    To prove Statement~\ref{item:GOLARecons}, note first that $M_{\textrm{g}} \geq \operatorname{len}(\phi_{0,n,r})$ for each $\phi_{0,n,r} \in \mathscr{F}(I^{\textrm{g}})$, in accordance with~\eqref{eq:MGlobal}.  Recalling that $X[m,n,r] = I^{\textrm{g}}[m,n,r]\langle x, \phi_{m,n,r} \rangle$ by~\eqref{eq:adaptInnerProduct}, we observe that the innermost summation of~\eqref{eq:golaReconstruction} is recognizable as an inverse discrete Fourier transform on $\mathbb{C}^{M_{\textrm{g}}}$, with $M_{\textrm{g}} = L / b_{\textrm{g}}$.  For fixed $n$ and $t \!- \!na \in \mathbb{Z}_{M_{\textrm{g}}}$, this term evaluates to $I^{\textrm{g}}[n,r] \mathcal{T}_{na}w_r[t]x[t]$:
    \begin{equation*}
        \label{eq:IDFT}
        \frac{1}{M_{\textrm{g}}} \!\! \sum_{m=0}^{M_{\textrm{g}}-1} \! X[\textstyle{\frac{M_L}{M_{\textrm{g}}}}m,n,r] e^{2\pi i m b_{\textrm{g}} t / L}
        \!\!\!\ \!\!\!\ \!\!\!\ \!\!\!\ \!\!\!\ \!\!\!\ \!\!\!\ \!\!
        \underset{\stackrel{~}{{\quad \footnotesize t-na \in \mathbb{Z}_{M_{\textrm{g}}}}}}{=}
        \!\!\!\ \!\!\!\ \!\!\!\ \!\!\!\ \!\!\!\ \!\!\!\ \!\!
        I^{\textrm{g}}[n,r] \mathcal{T}_{na}w_r[t]x[t] \text{.}
    \end{equation*}
    Perfect reconstruction then follows from the generalized overlap-add constraint of~\eqref{eq:adaptiveOlaConstraint}, as
    \begin{equation*}
        \frac{a}{\widehat{w}[0]} \sum_{n=0}^{N-1}\sum_{r=0}^{N-1} I^{\textrm{g}}[n,r] \mathcal{T}_{na}w_r[t]x[t]
        = \frac{a}{\widehat{w}[0]} \frac{\widehat{w}[0]}{a} x[t]  = x[t] \text{.}
     \end{equation*}

    For the case of a local selection function $I^l$, note that the generalized overlap-add constraint of~\eqref{eq:adaptiveOlaConstraint} is still implied by~\eqref{eq:olaConstraint}, since the argument for the case of admissible $I^{\textrm{g}}$ holds independently of the modulation structure employed. Consequently, by substituting $M_r, b_r$ for $M_{\textrm{g}}, b_{\textrm{g}}$ and noting that $M_r \geq \operatorname{len}(\phi_{0,n,r})$ by~\eqref{eq:MLocal} for each $\phi_{0,n,r} \in \mathscr{F}(I^l)$, we see that the result of~\eqref{eq:golaReconstruction} also holds for all $\mathscr{F}(I^l)$.
\end{IEEEproof}

\subsection{Reconstruction via Canonical Dual Superposition Frames}
\label{sec:canonicalDual}

We next develop the reconstruction properties of our superposition families in a frame-theoretic context, noting that they qualify as ``painless nonorthogonal expansions''~\cite{Daubechies86}, and that the development below also holds for more general nonstationary Gabor frames~\cite{JailletDorfler09}.  In analogy to Definition~\ref{def:frameOpr}, we associate a superposition frame operator $S_I: \mathbb{C}^L \rightarrow \mathbb{C}^L$ through its action on any $x \in \mathbb{C}^L$ as $S_I x = \sum_{\phi \in \mathscr{F}(I)} \langle x, \phi_{m,n,r} \rangle \phi_{m,n,r}$.
\begin{defn}[Superposition$\!$ Frame$\!$  Operator,$\!$  Walnut$\!$  Form]\label{def:superposWalnut}
The discrete Walnut representations of superposition frame operators $S_{I^{\textrm{g}}}$ and $S_{I^l}$ are respectively given by the $L \times L$ positive semi-definite matrices with entries
\begin{align*}
    S_{I^{\textrm{g}}}[t,t'] & \triangleq \! M_{\textrm{g}} \, \mathbb{I}_{M_{\textrm{g}}\backslash(t-t')}[t\!-\!t']
    \!\! \sum_{n,r=0}^{N-1} I^{\textrm{g}}[n,r] \mathcal{T}_{na}w_r[t] \overline{\mathcal{T}_{na}w_r[t']}
    \text{,} \\
    S_{I^l}[t,t'] & \triangleq \! \sum_{n,r=0}^{N-1} I^l[n,r] M_r \mathbb{I}_{M_r\backslash(t-t')}[t\!-\!t'] \mathcal{T}_{na}w_r[t]\overline{\mathcal{T}_{na}w_r[t']} \text{.}
    \end{align*}
\end{defn}

Theorem~\ref{thm:adaptiveFrames} implies that any Gabor frame $\mathscr{G}(w,a,b)$ and admissible $I^{\textrm{g}}$ or $I^l$ together give rise to a superposition frame, and hence the corresponding superposition frame operators are of full rank.  Thus, to each global superposition frame $\mathscr{F}(I^{\textrm{g}}) = \{ \phi_{m,n,r} \}$ corresponds a unique canonical dual frame $\{ \widetilde{\phi}_{m,n,r} \}$, whose elements are obtained in turn as $\widetilde{\phi}_{m,n,r} \triangleq S_{I^{\textrm{g}}}^{-1} \phi_{m,n,r}$.  Accordingly, when $S_{I^{\textrm{g}}}$ is diagonal we may index elements of $\{ \widetilde{\phi}_{m,n,r} \}$ by the same admissible $I^{\textrm{g}}$, and we obtain the following reconstruction property: %
\begin{equation*}
\forall x \in \mathbb{C}^L, t \in \mathbb{Z}_L, \,\,\,
    x[t] =
    \!\!\!\! \!\!\!\! \!\!\!\! \!\!\!\!
    \sum_{\qquad m,n,r:I^{\textrm{g}}[m,n,r]=1}
    \!\!\!\! \!\!\!\! \!\!\!\!
    \langle x, \phi_{m,n,r} \rangle \widetilde{\phi}_{m,n,r}[t]
     \text{,}
\end{equation*}
with the above also holding for local $I^l$ by Theorem~\ref{thm:adaptiveFrames}. We thus denote by $\widetilde{\mathscr{F}(I^{\textrm{g}})}$ or $\widetilde{\mathscr{F}(I^l)}$ the corresponding dual frames, and observe the following consequence of the Walnut representation of Definition~\ref{def:superposWalnut} above.

\begin{thm}[Fast Inversion via Canonical Dual]\label{prop:diagDual}
For any $\mathscr{F}(I^{\textrm{g}})$, $\mathscr{F}(I^l)$ derived from a Gabor frame $\mathscr{G}(w,a,b)$ on $\mathbb{C}^L$, the corresponding operators $S_{I^{\textrm{g}}}$ and $S_{I^l}$ are \emph{diagonal}, and each canonical dual frame element takes the form
\begin{equation}
    \label{eq:dualWalnutFormula}
    \widetilde{\phi}_{m,n,r}[t] = \frac{\mathcal{M}_{mb_L}\mathcal{T}_{na}w_r[t]}{M_{\textrm{g}} \textstyle{\sum_{n'=0}^{N-1}\sum_{r'=0}^{N-1}} I^{\textrm{g}}[n',r']\left|\mathcal{T}_{n'a}w_{r'}[t]\right|^2}
\end{equation}
for $\mathscr{F}(I^{\textrm{g}})$, and similarly for $\mathscr{F}(I^l)$ with respect to each $M_r$.
\end{thm}

Note that the corresponding formula for the nonstationary Gabor frames of~\cite{JailletDorfler09} in the diagonal case is similar to~\eqref{eq:dualWalnutFormula}. However, in the superposition frame setting, the constraints on the window structure not only preserve lower frame bounds and yield fast inversion via the constant overlap-add method, but also enable signal-independent  evaluation of the canonical dual in certain cases, as we now show.

\subsection{Adaptive Lapped Superposition Frames}
\label{sec:canonicalDualFurther}

Reconstruction via the canonical dual $\widetilde{\mathscr{F}(I^{\textrm{g}})}$ according to~\eqref{eq:dualWalnutFormula} requires knowledge of the admissible selection function $I^{\textrm{g}}[n,r]$ corresponding to a given signal adaptation.  This stands in contrast not only to the usual Gabor setting, wherein the form of the canonical dual frame can be obtained immediately, but also to the constant overlap-add approach described in Section~\ref{sec:Ola}, which avoids computation of the canonical dual entirely.  However, by coupling our superposition construction with the following \emph{neighbor overlap} condition, we are able to compute $\widetilde{\mathscr{F}(I^{\textrm{g}})}$ \emph{prior} to adaptation---that is, without knowledge of which ordered partition function will be used in subsequent signal analysis.

\begin{defn}[Neighbor$\!$ Overlap$\!$ Condition]\label{defn:limitOverlap}
    Let $\mathscr{G}(w, a, \cdot)$ be a Gabor system on $\mathbb{C}^L$, with $N = L/a$.  It is said to satisfy the \emph{neighbor overlap} condition if, for all $n,n'\in\mathbb{Z}_N$,
    \begin{equation}
     \label{eq:limitOverlap}
     \operatorname{supp} \left ( \mathcal{T}_{na}w \right ) \, \cap \, \operatorname{supp} \left ( \mathcal{T}_{n'a}w \right ) = \emptyset
     \,\, \text{if} \,\, |n-n'| > 1
     \text{.}
    \end{equation}
\end{defn}
Any admissible selection function preserves the neighbor overlap property, leading to the following notion of \emph{lapped superposition frames}, whose properties we develop below.
\begin{defn}[Adaptive Lapped Superposition Frames]
        Let $\mathscr{G}(w, a, \cdot)$ be a Gabor frame on $\mathbb{C}^L$ that \emph{simultaneously} satisfies the overlap-add constraint of~\eqref{eq:olaConstraint} and the neighbor-overlap condition of~\eqref{eq:limitOverlap}. Then for any admissible $I^{\textrm{g}}$, we call $\mathscr{F}(I^{\textrm{g}})$ an adaptive \emph{lapped superposition frame}.
\end{defn}

Note that the overlap-add constraint of~\eqref{eq:olaConstraint} ensures a partition of unity by window translates, while the neighbor overlap condition of~\eqref{eq:limitOverlap} is also required in the case of lapped \emph{orthogonal} transforms (see, e.g.,~\cite{Mallat99}).  While our construction retains the flavor of time-varying lapped transforms~\cite{herley1993ttf, Wang06}, we emphasize that the resultant frames can avoid the lack of translation invariance inherent in the orthogonal setting, while still ensuring fast reconstruction.

\begin{figure}[!t]
  \centering
  \includegraphics[width=.74\columnwidth]{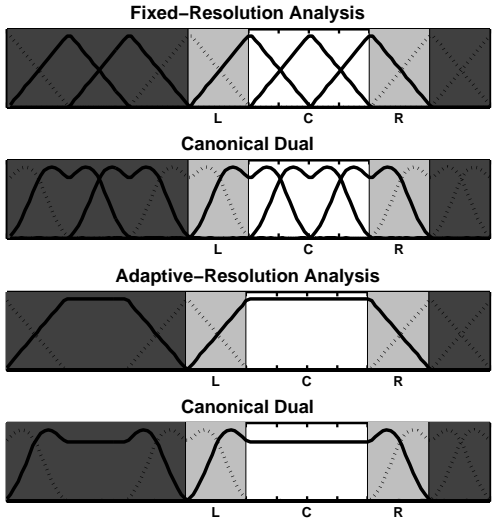}
  \caption{\label{fig:olaCanonicalDual}Translates $\{\mathcal{T}_{na}w\}$ from a Gabor frame $\mathscr{G}(w,a,\cdot)$ constructed from triangular windows with $50\%$ overlap (top panel), shown with translates of its canonical dual window $\widetilde{w}$ (second panel); remaining panels repeat this sequence for an induced superposition frame $\mathscr{F}(I^{\textrm{g}})$. On the left and right sets $L$ and $R$, the canonical dual windows of $\widetilde{\mathscr{F}(I^{\textrm{g}})}$ agree pointwise with those of $\mathscr{G}(\widetilde{w},a,\cdot)$; on the center set $C$, that of $\widetilde{\mathscr{F}(I^{\textrm{g}})}$ is constant.}
\end{figure}

We show below that if $\mathscr{F}(I^{\textrm{g}})$ is a lapped superposition frame derived from a Gabor frame $\mathscr{G}(w, a, b)$, then its canonical dual frame elements may be pre-computed.  This situation is illustrated in Fig.~\ref{fig:olaCanonicalDual},
where the support sets of a window $w$, its canonical dual $\widetilde{w}$, and their immediate neighbors are partitioned into subsets labeled $L$, $C$, and $R$.  Since $\mathscr{F}(I^{\textrm{g}})$ must inherit the neighbor-overlap condition from $\mathscr{G}(w, a, b)$, it follows that whenever $\mathscr{F}(I^{\textrm{g}})$ admits a diagonal frame operator, the corresponding canonical dual windows of $\widetilde{\mathscr{F}(I^{\textrm{g}})}$ are \emph{constant} on the center set $C$, as shown in Fig.~\ref{fig:olaCanonicalDual}.  Moreover, the highlighted dual superposition window is pointwise equal to $\widetilde{w}$ on the sets $L$ and $R$ (see bottom two panels of Fig.~\ref{fig:olaCanonicalDual}).

To formalize this intuition, define for each $r$ the sets
\begin{equation}\label{eqn:sets}
    \begin{cases}
    L_r \triangleq \operatorname{supp}(w_r) \cap \operatorname{supp}(\mathcal{T}_{-(r+1)a} w_r) \text{,} & \\
    R_r \triangleq \operatorname{supp}(w_r) \cap \operatorname{supp}(\mathcal{T}_{(r+1)a} w_r) \text{,} & \\
    C_r \triangleq \operatorname{supp}(w_r) \setminus (L_r \cup R_r) \text{,}
    \end{cases}
\end{equation}
and note that, for any global selection function $I^{\textrm{g}}$, $\phi_{\cdot,n,r} \in \mathscr{F}(I^{\textrm{g}})$ and its dual $S_{I^{\textrm{g}}}^{-1}\phi_{\cdot,n,r} \in \widetilde{\mathscr{F}(I^{\textrm{g}})}$ are \emph{both} supported exclusively on the set $\mathcal{T}_{na}(L_r \cup C_r \cup R_r)$. Then the following theorem, proved in the appendix, establishes our main result: the canonical dual frame of any $\mathscr{F}(I^{\textrm{g}})$ can be computed independently of any ordered partition function, and, therefore, can be computed \emph{prior} to observing any data.

\begin{thm}[Canonical Duals of Adaptive Lapped Frames]
    \label{thm:golaStructure}
    Let $\mathscr{F}(I^{\textrm{g}})$ arise from a Gabor frame $\mathscr{G}(w, a, \cdot)$ for $\mathbb{C}^L$
    satisfying~\eqref{eq:olaConstraint} and~\eqref{eq:limitOverlap}.
    Then every $\widetilde{\phi}_{m,n,r} \in \widetilde{\mathscr{F}(I^{\textrm{g}})}$ can be constructed by modulations $\mathcal{M}_{mb_L}$ and translations $\mathcal{T}_{na}$ of lapped windows corresponding to each $r$ as follows:
        \begin{equation*}
        \widetilde{\phi}_{0,0,r}[t] \triangleq\, \frac{1}{M_{\textrm{g}}} \textstyle \!\cdot\!
          \begin{cases}
            \frac{w[t]}{\sum_{n=0}^{N-1}|w[t-na]|^2} & \text{if $t \in L_r$,} \\
            \frac{a}{\widehat{w}[0]} & \text{if $t \in C_r$,} \\
            \frac{w[t-(r+1)a]}{\sum_{n=0}^{N-1}|w[t-na]|^2} & \text{if $t \in R_r$.}
          \end{cases}
        \end{equation*}
    Here the sets $L_r,C_r,R_r$ are defined in~\eqref{eqn:sets}, and we note that only the expression for $C_r$ is to be employed when $\mathscr{F}(I^{\textrm{g}})$ is comprised entirely of modulations of $w_{N-1}$.
\end{thm}

We note that many popular Gabor systems $\mathscr{G}(w, a, \cdot)$ satisfy the requirements of this theorem---including triangular, Hamming, and raised-cosine windows $w$ at $50\%$ overlap---thus enabling a variety of new adaptive, lapped superposition frame families that all admit fast reconstruction.

\begin{table*}[t]
\caption{\label{tb:complexityTable}Computational complexity orders $\mathcal{O}(\cdot)$ of various analysis and synthesis algorithms considered in the article}
\centering
\begin{tabular}{rc|cc|c|c}
    & \multicolumn{2}{c}{\underline{No Adaptation}} & \multicolumn{3}{c}{\underline{Adaptation}} \\
    Analysis Complexity  & \multicolumn{2}{c}{$NM(1+\log_2M)$}  & \multicolumn{3}{c}{$N_\textrm{g}M_\textrm{g}(1+\log_2M_\textrm{g})$} \\
    Synthesis Method    & Overlap-Add & Canonical Dual &  Section~\ref{sec:Ola} & Section~\ref{sec:canonicalDual} & Section~\ref{sec:canonicalDualFurther} \\
    Synthesis Complexity & $NM\log_2M$ & $NM(1+\log_2M)$  & $N_\textrm{g}M_\textrm{g}(\log_2M_\textrm{g})$ & $N_\mathrm{g}M_\mathrm{g
    }(3 + \log_2M_\textrm{g})$ & $N_\textrm{g}M_\textrm{g}(1+\log_2M_\textrm{g})$  \\
\end{tabular}
\vspace{-\baselineskip}%
\end{table*}

\subsection{Adaptive Dyadic Superposition Frames}
\label{sec:dyadicSuperpositionFrames}

Here we construct a class of so-called \emph{dyadic} superposition frames that admit offline canonical dual construction even when the overlap-add constraint of~\eqref{eq:olaConstraint} is \emph{not} satisfied.  We base this construction on the notion of \emph{dyadic} ordered partition functions, which may be thought of as indexing binary trees.

\begin{defn}[Dyadic Admissible Selection Functions]
   \label{exmp:binaryTree}
    Let the number of translates $N$ of $w$ in $\mathscr{G}(w,a,b)$ be a power of two, and define the set $\mathscr{H} \triangleq \{0,1,\ldots,\log_2 N \}$.
    An ordered partition function $\widetilde{I}^d[n,r]$ is \emph{dyadic} if it satisfies the conditions of Definition~\ref{def:stickBreaking} and
    \begin{equation*}
        \widetilde{I}^d[n,r] = 1 \,\, \text{\emph{only} if} \,\, r = 2^{h}-1 \,\, \text{for some} \,\, h \in \mathscr{H} \text{.}
    \end{equation*}
    We denote by $I^d[m,n,r]$ a \emph{dyadic} admissible selection function induced by $\widetilde{I}^{d}[n,r]$ and a global frequency lattice constant $b_\textrm{g} = L/M_\textrm{g}$, with $M_\textrm{g}$ defined according to~\eqref{eq:MGlobal}.
\end{defn}

Viewing the $N$ translates of $w$ in $\mathscr{G}(w,a,b)$ as leaves of a binary tree of height $\log_2N$, a dyadic ordered partition function selects windows corresponding to some tree level $h$.

\begin{defn}[Dyadic Gabor and Superposition Frames]\label{defn:dySupFr}
 Fix an initial Gabor frame $\mathscr{G}(w, a, b)$, such that $N = L/a$ is a power of two, and let $h \in \mathscr{H}$ index \emph{height} in a binary tree.  Now restrict $r \in \mathbb{Z}_n$ to the index set $\mathscr{R} \triangleq \{2^h-1\}$, and fix for each $r \in \mathscr{R}$ a time lattice constant $a_r \triangleq a(r+1)$, and a dyadic admissible selection function $I^d$ with associated global frequency lattice constant $b_\textrm{g}$. Then:
 \begin{enumerate}

    \item We define \emph{dyadic superposition Gabor frame} $G^d_r(b_\textrm{g})$ for all  $r \in \mathscr{R}$, and their union $G^d_{\cup r}(b_\textrm{g})$, as follows:
    \begin{align*}
        G^d_r(b_\textrm{g}) & \triangleq \mathscr{G}(w_r, a_r, b_\textrm{g}) \subseteq \mathscr{G}(w_r, a, b_L)
        \text{,}
        \\
        G^d_{\cup r}(b_\textrm{g}) & \triangleq \displaystyle{\cup_{r \in \mathscr{R}}} G_{r}^d(b_\textrm{g}) \subseteq \cup_{r=1}^{N-1}\mathscr{G}(w_r, a, b_L)
        \text{.}
    \end{align*}

    \item We call $\mathscr{F}(I^d) \!\subset\! G^d_{\cup r}(b_{\textrm{g}})$ a \emph{dyadic superposition frame}.
 \end{enumerate}
\end{defn}
The fact that $\mathscr{F}(I^d)$ and every dyadic $G^d_r(b_\textrm{g})  = \mathscr{G}(w_r, a_r, b_\textrm{g})$ are frames for $\mathbb{C}^L$ follows from the assumption that $\mathscr{G}(w, a, b)$ is a frame, coupled with the result of Theorem~\ref{thm:adaptiveFrames}.

An example of this construction is illustrated in Fig.~\ref{fig:dyadicExample}:
\begin{figure}[!t]
  \centering
  \includegraphics[width=\columnwidth]{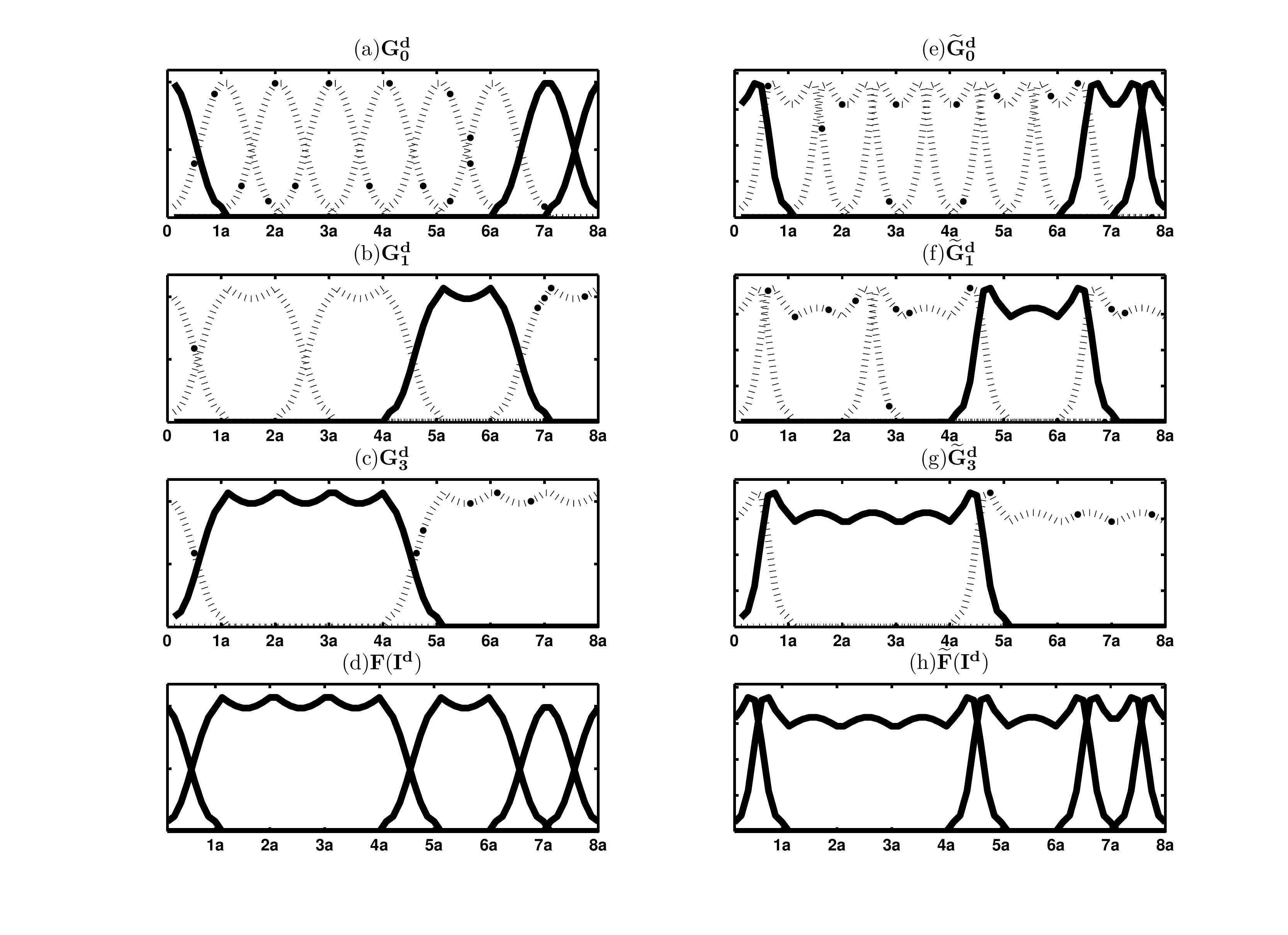}
  \vspace{-0.75\baselineskip}%
  \caption{\label{fig:dyadicExample}Repeated superposition merges of Hamming windows, following the structure of a binary tree, are shown in panels (a-c), along with the canonical dual windows associated to each corresponding Gabor system (e-g).  A dyadic superposition frame $\mathscr{F}(I^d)$ can be formed from the selected unmodulated elements of (d), which in turn will admit the corresponding canonical dual windows shown in (h), computed according to~\eqref{eq:dualWalnutFormula}.}
\end{figure}
The left panel shows Hamming windows $w \equiv w_0$ at 50\% overlap, along with the corresponding dyadically-indexed superposition windows $w_1$ and $w_3$, and an example dyadic superposition frame $\mathscr{F}(I^d)$ at bottom.  The right panel shows canonical duals associated to each dyadic superposition Gabor frame $G^d_r(\cdot)$, for $r \in \{0,1,3\}$, along with the corresponding canonical dual $\widetilde{\mathscr{F}(I^d)}$ at bottom right.  The fact that $\widetilde{\mathscr{F}(I^d)}$ contains elements from each of these individual dual Gabor frames $\widetilde{G^d_r(\cdot)}$ is verified by the following theorem.

\begin{thm}[Canonical Duals of Adaptive Dyadic Frames]
    \label{thm:dyadicStructure}
    Let a dyadic superposition frame $\mathscr{F}(I^d) \subset G^d_{\cup r}(b_{\textrm{g}})$ arise from a Gabor frame $\mathscr{G}(w, a, b)$ for $\mathbb{C}^L$ satisfying the neighbor-overlap condition of~\eqref{eq:limitOverlap}. Then its canonical dual $\widetilde{\mathscr{F}(I^d)}$ can be computed by either path of the following commutative diagram, where $S_{I^d}$ is the frame operator associated to $\mathscr{F}(I^d)$, and $S^d_r$ is that associated to each dyadic Gabor frame $G^d_r$:
    \begin{equation}
        \label{eq:cdDyadic}
            \begin{CD}
            G^d_{\cup r}(b_\textrm{g}) @> (S^d_r)^{-1}, r \in \mathscr{R} >> \cup_{r \in \mathscr{R}} \widetilde{G^d_r(b_\textrm{g})} \\ @VVI^dV @VVI^dV\\
             \mathscr{F}(I^d) @> (S_{I^d})^{-1}>> \widetilde{\mathscr{F}(I^d)}
            \end{CD}\text{.}
    \end{equation}
\end{thm}

Observe that computing the $\widetilde{\mathscr{F}(I^d)}$  via direct inversion of its (diagonal) frame operator $S_{I^d}$ corresponds to the down-and-right path in~\eqref{eq:cdDyadic}, and hence requires knowledge of $I^d$. However, Theorem~\ref{thm:dyadicStructure} implies that all elements in $\widetilde{\mathscr{F}(I^d)}$ can be pre-computed by instead following the right-and-down path. As in the case of lapped superposition frames in Section~\ref{sec:canonicalDualFurther}, the neighbor overlap condition of~\eqref{eq:limitOverlap} plays a key role.

\subsection{Computational Complexity}
\label{sec:compComplexity}

We now address the relative computational complexity of the various reconstruction methods presented above.  Recall that these algorithms all yield a diagonal frame operator, implying that only $\mathcal{O}(L)$ operations are required for its inversion, compared to $\mathcal{O}(L^3)$ in the worst case of a non-diagonal frame operator lacking any special structure.  While in some circumstances this complexity can be reduced (see, e.g.,~\cite{Qiu1995}), such schemes remain super-linear in $L$, rendering them impractical for use in applications with $L \gg 1$.

Recall that when no windows are merged, we recover a Gabor frame $\mathscr{G}(w,a,b)$ containing $NM$ elements, with $Na = Mb = L$. If the associated frame operator is diagonal, an FFT-based approach requires $NM(1 + \log_2M)$ complex multiplications to compute the short-time analysis coefficients $X[m,n]$ from $x \in \mathbb{C}^L$, with $NM$ of these needed to obtain the $N$ individual short-time segments $\{\mathcal{T}_{na}w[t]x[t]\mathbb{I}_{\operatorname{supp}(\mathcal{T}_{na}w)}[t] \}$.  When reconstruction proceeds via inversion of the frame operator, then $NM\log_2M$ operations are required to compute the necessary inverse FFTs, plus $NM$ operations to multiply each resultant segment by the appropriate dual window.  In the overlap-add setting, this window is the identity, and hence these latter $MN$ operations are avoided.
An in-depth discussion of the general (non-diagonal) case, including methods based on the Zak transform, is given in~\cite{Qiu1995} and~\cite{Orr1993}.

In the case that windows are merged, note that a global superposition frame $\mathscr{F}(I^{\textrm{g}})$ derived from $\mathscr{G}(w,a,b)$ comprises $M_\textrm{g} = L/b_{\textrm{g}} \geq M$ modulations of $N_\textrm{g} \leq N$ windows:
\begin{equation}
    \label{eq:numWindows}
    N_\textrm{g} = \sum_{n=0}^{N-1}\sum_{r=0}^{N-1} I^{\textrm{g}}[n,r] \text{.}
\end{equation}
It follows that $N_\textrm{g}M_\textrm{g}(\log_2 M_\textrm{g} + 1)$ complex multiplications are required to compute the analysis coefficients $X[m,n,r]$ via $N_\textrm{g}$ FFTs, followed by $N_\textrm{g}M_\textrm{g}\log_2M_\textrm{g}$ operations required for the $N_\textrm{g}$ inverse FFTs required for reconstruction.

If $\mathscr{F}(I^\textrm{g})$ arises from the special case described in Theorem~\ref{thm:golaStructure}, then the elements of $\widetilde{\mathscr{F}(I^\textrm{g})}$ can be pre-computed, and only $N_\textrm{g}M_\textrm{g}$ extra operations are necessary for multiplication by the requisite canonical dual windows. In the general case, elements of $\widetilde{\mathscr{F}(I^\textrm{g})}$ must be computed directly via~\eqref{eq:dualWalnutFormula}, as a function of the chosen selection function $I^\textrm{g}$.  Appealing to~\eqref{eq:numWindows}, we see that this computation can be accomplished using another $2N_\textrm{g}M_\textrm{g}$ calculations, leading to a total of $N_\textrm{g}M_\textrm{g}(3 + \log_2M_\textrm{g})$ complex multiplications. This analysis is also an upper bound for the worst-case complexity for any $\mathscr{F}(I^l)$ having the same ordered partition function as $\mathscr{F}(I^\textrm{g})$.

The various analysis and synthesis complexities discussed above are summarized in Table~\ref{tb:complexityTable}. In practice, the complexity of the signal adaptation procedure must also be taken into account. This complexity depends both on the method for searching among ordered partition functions $\widetilde{I}[n,r]$, and the cost function used to compare them. Since it is clearly infeasible to compare all $2^{N-1}$ ordered partition functions by exhaustive search, we next consider greedy and dynamic-programming-based approaches below.  In both cases, the complexity of evaluating the associated cost functions increases with window length, and therefore it is advisable in practice to set an upper bound on the maximal number of window merges.

\section{Signal Adaptation Algorithms and Examples}
\label{sec:algSpecifics}

Signal-adaptive modification of an initial Gabor frame $\mathscr{G}(w, a ,b)$ on $\mathbb{C}^L$ via superposition can produce any one of the possible $2^{N-1}$ superposition frames whose properties we characterized in Sections~\ref{sec:adaptiveFrameAnalysis} and~\ref{sec:synthesis}. We now detail two instances of a broad class of signal adaptation algorithms, any of which can be used to select a superposition frame for subsequent signal analysis. We propose both greedy and dynamic programming approaches in Section~\ref{sec:algorithms}, and illustrate their performance with two brief examples in Section~\ref{sec:Examples}.

\subsection{Signal-Adaptive Superposition Frame Selection}
\label{sec:algorithms}

The first of the two Gabor frame adaptation algorithms we describe is a simple greedy approach that can be implemented by ``growing'' a given window forward in time through successive attempts to merge it with its subsequent neighboring translates~\cite{Rudoy2008}. Whenever a proposed merge fails, the procedure resets and repeats, halting when the end of the data stream is reached (or, equivalently in our cyclic setting, when the initial window is once again encountered).

A decision whether or not to merge adjacent windows can be made based on any suitable cost function. As one example, we employ the time-frequency concentration measure appearing in the popular work of~\cite{JonesBaraniuk94, Jones1995} on adaptive optimal-kernel time-frequency representations. Specifically, consider a short-time segment $x_{n,r}[t] \triangleq \{\mathcal{T}_{na}w[t]x[t]\mathbb{I}_{\operatorname{supp}(\mathcal{T}_{na}w)}[t] \}$, and define its time-frequency concentration in the manner of~\cite{Jones1990, JonesBaraniuk94, Jones1995, Rudoy2008}
\begin{equation}
\label{eq:kurtosis}
    C(x_{n,r}) \triangleq \frac{\sum_{m=0}^{M_r-1} \left| \langle x, \mathcal{M}_{mb_r}\mathcal{T}_{na}w_r \rangle\right|^4}
    {\left(\sum_{m=0}^{M_r-1} \left| \langle x, \mathcal{M}_{mb_r}\mathcal{T}_{na}w_r \rangle \right|^2\right)^2},
\end{equation}
with $M_r,b_r$ defined via~\eqref{eq:MLocal}.
This ratio of powers of norms of short-time Fourier coefficients is suggestive of an ``empirical spectral kurtosis,'' and has also been used in minimum entropy deconvolution~\cite{Wiggins78}; other choices are also possible~\cite{Nesbit2009}.

As shown in~\cite{JonesBaraniuk94}, maximizing~\eqref{eq:kurtosis} favors short-time segments that \emph{concentrate} local signal energy within the smallest regions of the time-frequency plane. Indeed, below we obtain similar results on an example akin to the one employed in~\cite{JonesBaraniuk94}:  the resultant superposition frames comprise shorter windows near time-localized transients, and longer windows near oscillatory signal portions. The resulting procedure requires $\mathcal{O}(N)$ iterations and is summarized in Algorithm~\ref{alg:kurtosis}.

\begin{algorithm}
{\small \caption{\label{alg:kurtosis} Adaptation via Greedy Selection~\cite{Rudoy2008}}
        \vspace{.2\baselineskip}%
        \textit{Initialization}
        \begin{itemize}
            \item Fix input data $x \in \mathbb{C}^L$ and a Gabor frame $\mathscr{G}(w,a,b)$

            \item Set $(p, n_p) = (0,0)$ and initialize $\widetilde{I}[n,r]$ to be the $N$-part ordered partition function of Example~\ref{ex:uniform}
        \end{itemize}
        \vspace{.1\baselineskip}%
        \textit{Greedy Selection:} For $n=0, 1, \ldots, N-1$,

        \begin{itemize}
            \item Compute a merged window
            \begin{equation*}
                \mathcal{T}_{n_p a}w_{p+1} = \mathcal{T}_{n_p a} w_p + \mathcal{T}_{na}w
            \end{equation*}
            and $C(x_{n_p,p+1})$, $C(x_{n_p,p})$ and $C(x_{n,0})$ via~\eqref{eq:kurtosis}

            \item If $C(x_{n_p,p+1}) \leq  \max(C(x_{n_p,p}), C(x_{n,0}))$,
            reject the proposed merge: set $(p, n_p)$ as $(p, n+p+1)$, and leave $\widetilde{I}[n,r]$ unchanged

            \item Otherwise accept the proposed merge: set $(p, n_p)$ as $(p+1, n_p)$ and update $\widetilde{I}[n,r]$ as
             \begin{align*}
                &\widetilde{I}[n_p, p+1] = 1 \quad \quad \quad \,\, \text{(add: $\mathcal{T}_{n_p a}w_{p+1}$),} \\
                &\widetilde{I}[n_p, p] = \widetilde{I}[n, 0] = 0 \quad \text{(remove: $\mathcal{T}_{n_p a} w_p, \, \mathcal{T}_{na}w$)}
            \end{align*}
        \end{itemize}
        \textit{Output:} Return the set of variable-length windows induced by $\widetilde{I}[n,r]$
    }%
\end{algorithm}

The second algorithm we present is based on the dynamic programming approach to adaptive segmentation popular in the audio coding literature~\cite{RamchandranOrchard97,PrandoniVetterli00,NiamutHeusdens05,HeusdensJensen2005, RodbroHeusdens06}. The basic idea is to fix an \emph{additive} cost function $J(\cdot)$, and find an optimal ordered partition function $\widetilde{I}^*[n,r]$ in the sense that it minimizes the sum of individual segment costs $J(x_{n,r}[t])$:
\begin{equation*}
    \widetilde{I}^*[n,r] \triangleq \underset{\widetilde{I}[n,r]}{\operatorname{arg min}} \!\!\! \sum_{n, r: \widetilde{I}[n,r]=1} \!\!\!\!\! \widetilde{I}[n,r] J(x_{n,r}[t]) \text{.}
\end{equation*}
Many choices for $J(\cdot)$ are possible, including rate-distortion cost functions~\cite{PrandoniVetterli00}, sparsity-inducing measures~\cite{Nesbit2009}, and the well-known entropy cost of~\cite{Coifman1992}, which we employ below.

To formalize our approach, define $J^*_{n}$ as the minimum cost among ordered partition functions on $\{0, 1, \ldots, na-1\}$, and let $J_{n,r} \triangleq J(x_{n,r}[t])$ represent the cost associated to covering the region $\{na, na+1, \ldots, (n+r)a\}$. The resulting dynamic program requires $\mathcal{O}(N^2)$ iterations and is summarized in Algorithm~\ref{alg:dp}.

\begin{algorithm}
{\small \caption{\label{alg:dp} Adaptation via Dynamic Programming~\cite{RamchandranOrchard97}}
        \vspace{.2\baselineskip}%
        \textit{Initialization}
        \begin{itemize}
            \item Fix input data $x \in \mathbb{C}^L$, a Gabor frame $\mathscr{G}(w,a,b)$ and initialize the cost function $J^*_0 = 0$
            \item For each $\mathcal{T}_{na}w \in \{T_{na}w: n \in \mathbb{Z}_N \}$ calculate
            the support set
            \[
                \mathscr{D}_{n} \triangleq \{ t : \mathcal{T}_{na}w[t] > \mathcal{T}_{n'a}w[t], n \neq n' \in \mathbb{Z}_N  \}
            \]
        \end{itemize}

        \textit{Dynamic Program}
        \begin{itemize}
            \item    For $n = 0, 1, \ldots, N-1$, compute sequentially the nth segmental cost and associated boundary by
        \begin{align*}
            J_n^* = \min_{0 \leq r < n} \left (J_r^* + J_{n,r} \right ) \qquad
            b_n^* = \operatorname{arg} \min_{0 \leq r < n} \left (J_r^* + J_{n,r} \right ) \text{,}
        \end{align*}
        with $J_{n,r}$ calculated using signal data supported on $\cup_{k = n}^{n+r} \mathscr{D}_k$
        \vspace{.2\baselineskip}%
        \item Compute the optimal selection function $\widetilde{I}^*[n,r]$ using $\{b_n^* : n \in \mathbb{Z}_N \}$ via the standard ``backtracking'' procedure~\cite{RamchandranOrchard97}
        \end{itemize}
        \textit{Output:} Return the set of variable-length windows induced by $\widetilde{I}^*[n,r]$
    }%
\end{algorithm}

Note that to preserve cost additivity in the presence of overlapping, non-orthogonal windows, Algorithm~\ref{alg:dp} evaluates $J(x_{n,r}[t])$ on regions \emph{smaller} than those covered by the corresponding windows, in a manner which recovers the approach of~\cite{PrandoniVetterli00} in the block-Fourier case.

Once a set of variable-length windows is obtained via any selection procedure returning an ordered partition function, a local or global modulation structure can be chosen via~\eqref{eq:MLocal} or~\eqref{eq:MGlobal}, respectively, in order to obtain a signal-adaptive superposition frame.  In practice, application-specific considerations are likely to play a role in superposition frame selection, and to this end we note that a variety of other algorithms and approaches are possible (see, e.g.,~\cite{Heusdens2006, Basu2009, RudoyQuatieriWolfe09}).

\begin{figure*}[t]
  \centering
  \includegraphics[width=2\columnwidth]{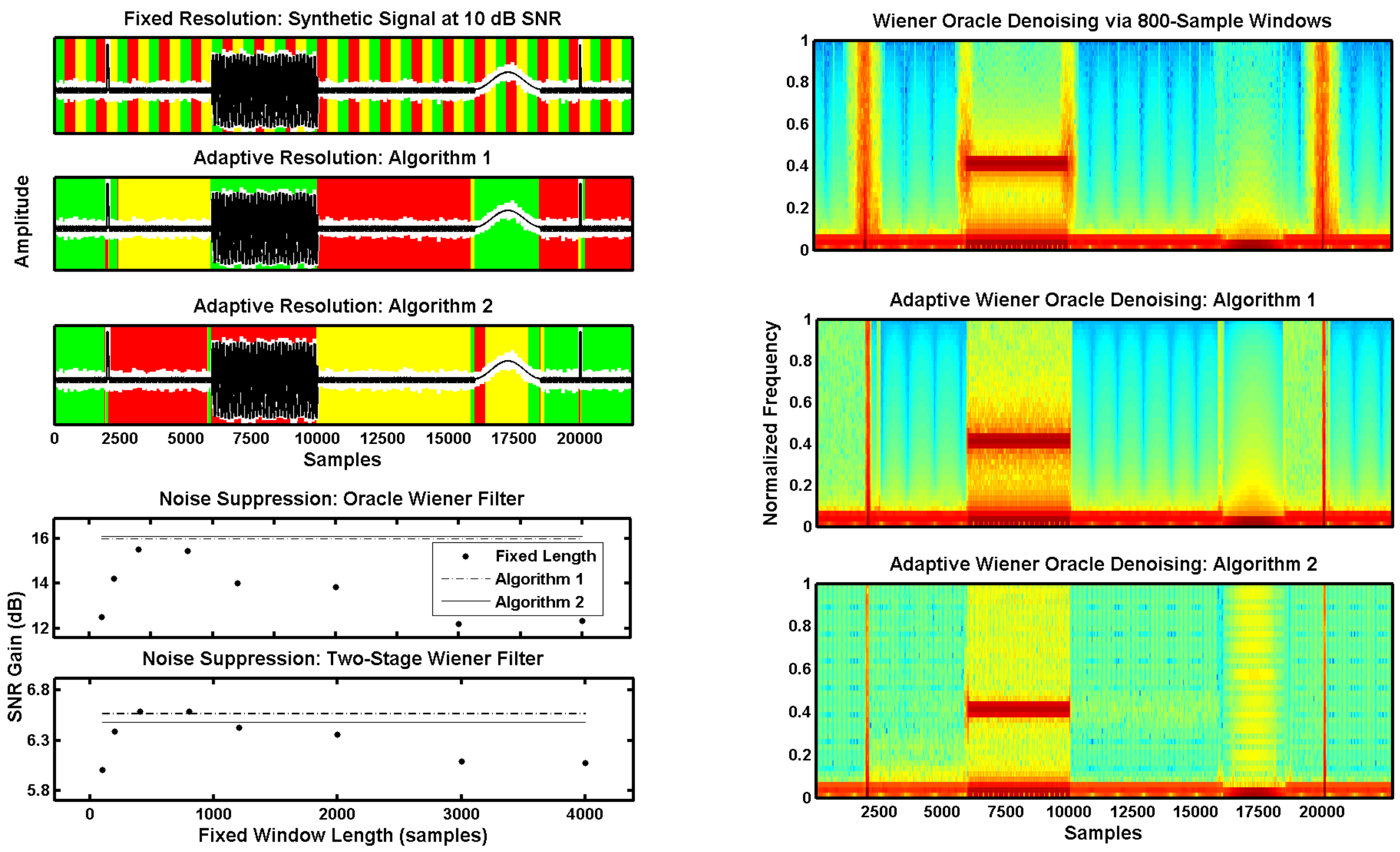}
  \vspace{-0.75\baselineskip}%
  \caption{\label{fig:synthPanel} Adaptive analysis-synthesis of a noisy synthetic signal ($10$~dB SNR) via superposition frames.  Top left: Adaptive segmentations of the \emph{noisy} signal (white, with clean version superimposed in black); background rectangles highlight the temporal extent of selected superposition windows.  Bottom left:  Results of a noise suppression experiment using an oracle (resp. two-stage) Wiener filter, averaged over $50$ trials.  Standard deviations range from $0.13$--$0.27$ dB (oracle) and from $0.03$--$0.06$ dB (two-stage).  Right: spectrograms formed from the oracle-denoised signals using $80$~sample Hanning windows with $50\%$ overlap.}
\end{figure*}

\subsection{Illustrative Examples}
\label{sec:Examples}

To conclude our investigation of superposition frames, we now consider two illustrative examples that combine Algorithms~\ref{alg:kurtosis} and~\ref{alg:dp} with the analysis and reconstruction procedures presented earlier.  These examples---a stylized synthetic waveform akin to the example employed in~\cite{JonesBaraniuk94} and a phonetically balanced speech utterance from the TIMIT corpus~\cite{Zue90}---both exhibit varying time-frequency structure, which in turn motivates signal-adaptive analysis and reconstruction.

Our first example signal $x$ (Fig.~\ref{fig:synthPanel}, top left) comprises a local and global sinusoidal term, two impulses, and a bump function.  We conducted a variety of experiments in which varying levels of white Gaussian noise $n$ were added to $x$, and Algorithms~\ref{alg:kurtosis} and~\ref{alg:dp} were then applied to $y \triangleq x + n$ to obtain signal-adaptive frame analysis coefficients $Y[m,n,r]$ on a global frequency lattice ($M_\textrm{g} = 6000$).  Using these as well as fixed-resolution analyses for a range of window lengths, we then applied to $Y[m,n,r]$ both an ``oracle'' Wiener suppression rule (local signal spectrum estimated by $|X[m,n,r]|^2$) and a two-stage Wiener suppression rule (by appropriately soft-thresholding $|Y[m,n,r]|^2$), and obtained a time-domain reconstruction $\widehat{x}$ via the corresponding canonical dual superposition frame.

The remainder of Fig.~\ref{fig:synthPanel} reports the results of a typical run at $10$~dB signal-to-noise ratio (SNR), with the superposition system of Algorithm~\ref{alg:kurtosis} derived from an initial Gabor system $\mathscr{G}(w,a,\cdot)$ comprising a $100$-sample Hamming window $w$ and time lattice constant $a=50$, and that of Algorithm~\ref{alg:dp} based on a $65$-sample Hamming window with $a=32$.  Although we have observed Algorithm~\ref{alg:dp} to be more noise-robust in practice, and to yield better performance with somewhat shorter initial windows, it may be seen that both algorithms yield broadly similar analyses with respect to dominant signal features at $10$~dB SNR.  Moreover, over a range of noise levels and fixed-resolution analyses, we have observed improved SNR gains $20\log_{10}(\|y-x\|/\|\widehat{x}-x\|)$ in both the oracle and two-stage cases, as shown in the bottom-left panels of Fig.~\ref{fig:synthPanel}.

Reconstruction spectrograms $20\log_{10} |\widehat{X}[m,n,r]|$---based on the \emph{oracle} denoising for visual clarity---are shown in the right-hand panel of Fig.~\ref{fig:synthPanel}.  They indicate that, in comparison to an \emph{a priori} well chosen fixed-resolution analysis using $800$-sample Hamming windows with $a = 400$, the onsets and offsets of localized time-frequency features are better preserved by superposition frames.  Since the best fixed-resolution window length is not known \emph{a priori} in practice, the adaptive approach remains attractive, despite the lessening SNR gains obtained in the simple two-stage denoising approach.  These results suggest the investigation of more sophisticated denoising schemes, also bearing in mind that in the case of nonstationary noise, the best adaptive analysis may well be SNR-dependent.

\begin{figure}
    \centering
    \includegraphics[width=\columnwidth]{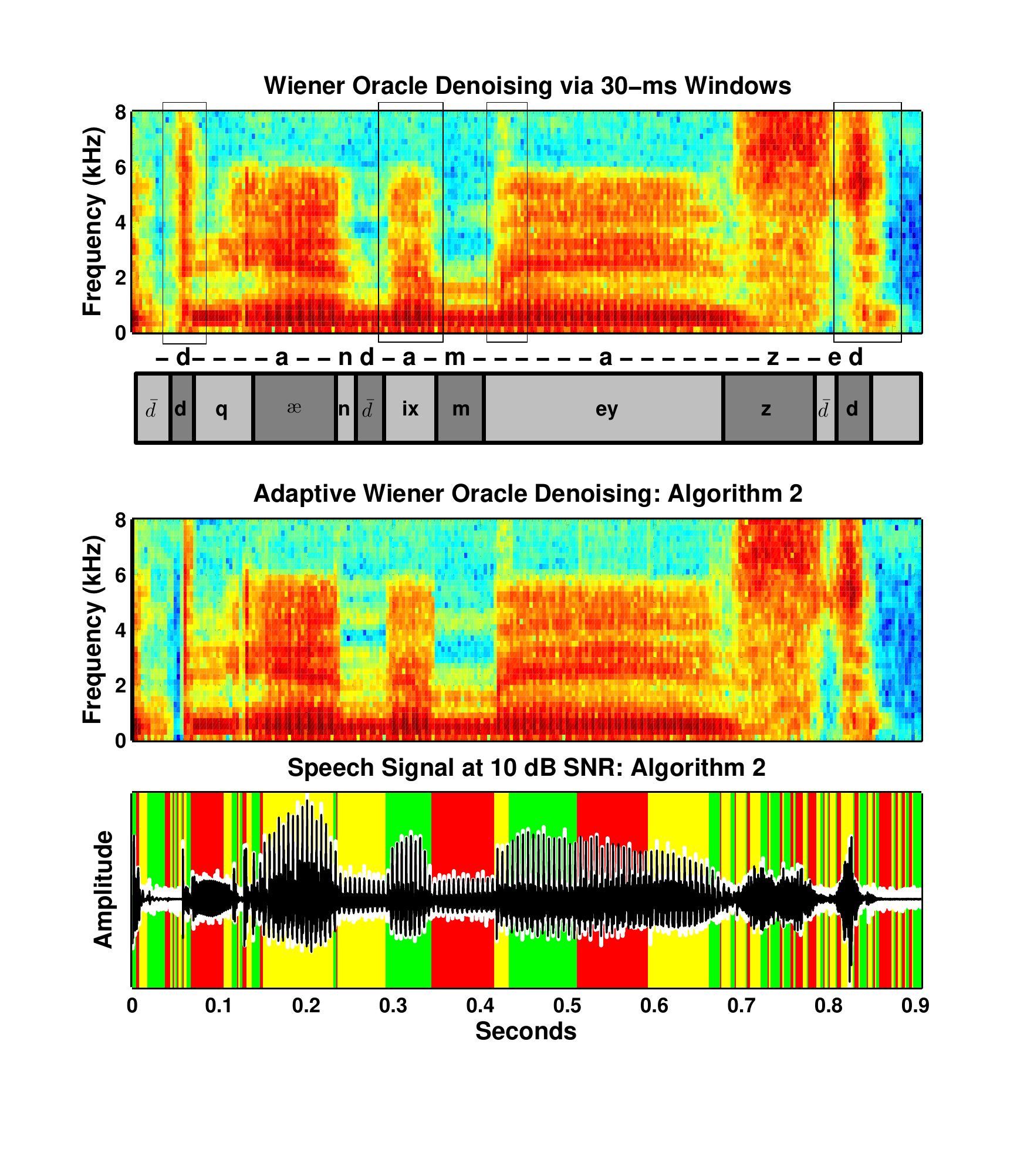}
    \vspace{-0.75\baselineskip}%
    \caption{\label{fig:speechPanel} Adaptive analysis-synthesis of a noisy speech signal ($10$~dB SNR).  Bottom: adaptive segmentation via Algorithm~\ref{alg:dp}, shown with spectrogram of adaptive oracle-Wiener-denoised version.  Top: spectrogram of fixed-resolution ($30$-ms) denoised version, shown with orthographic and phonetic transcriptions; boxes highlight temporal smearing of plosives and vowel onsets.  Both spectrograms were formed using $5$~ms Hanning windows with $50\%$ overlap.
    }
\end{figure}

We repeated the same battery of tests with our second example signal $x$ (Fig.~\ref{fig:speechPanel}, bottom), a portion of the phonetically-balanced TIMIT speech waveform /train/dr1/fsah0/si1244.wav corresponding to the phrase ``\ldots [eye]d and amazed.'' This utterance, chosen to illustrate time-varying spectral content typical of speech, contains two plosives ([eye]\emph{d}, amaze\emph{d}), two steady vowels (\emph{a}nd, \emph{a}mazed), and a time-varying diphthong (am\emph{a}zed).  The exact phonetic TIMIT segmentation (si1244.phn) is shown between the two spectrogram panels of Fig.~\ref{fig:speechPanel}, which correspond respectively to reconstructions based on an \emph{a priori} well chosen fixed-resolution (top, $30$~ms) and adaptive-resolution (middle, Algorithm~\ref{alg:dp}, starting from $3$~ms Hamming windows with $a = 1.5$~ms) oracle-Wiener denoising, respectively, at $10$~dB SNR.  The bottom panel of Fig.~\ref{fig:speechPanel} illustrates the corresponding adaptive analysis, which is seen to agree well with major features of the given TIMIT segmentation; Algorithm~\ref{alg:kurtosis} also yielded a similar analysis.

Following the same experimental procedure as in the case of Fig.~\ref{fig:synthPanel}, we observed broadly similar results---though with lower overall SNR gains obtained for this less stylized example.  Importantly, however, superposition windows are seen to better preserve vowel onsets and plosives; see in particular the boxed regions of the fixed-resolution spectrogram in the top panel of Fig.~\ref{fig:speechPanel}, corresponding to the two plosives and initial vowel-diphthong onsets in the word ``amazed.''

In this manner we see that superposition frames, when coupled with appropriate waveform adaptation criteria, show strong potential for use in a variety of signal-adaptive analysis-synthesis settings.  For signal enhancement applications, a natural next step would be to extend the approach of~\cite{Heusdens2006}, in which adaptive segmentation is used to estimate the local signal spectrum for enhancement purposes, but reconstruction is done using a fixed-resolution time-frequency lattice.  A variety of other multi-stage or iterative approaches suggest themselves, given the additional flexibility engendered by the overcomplete, signal-adaptive superposition frames presented in this article.

\section{Discussion}
\label{sec:discussion}

In this article we have introduced a broad family of adaptive, linear time-frequency representations termed superposition frames, and showed that they admit a host of desirable properties, including fast overlap-add reconstruction akin to standard short-time Fourier techniques.  Through a discussion of signal adaptation criteria and multiple examples, the resultant analysis-synthesis systems were seen to provide an effective and practical method for realizing signal-adaptive time-frequency analysis coupled with fast reconstruction.

Relative to other adaptive time-frequency methods, a number of open questions remain.  First, while many aspects of our construction admit straightforward extension to other Hilbert spaces of interest such as $\ell^2(\mathbb{Z})$ or $L^2(\mathbb{R})$, the present article has not addressed window requirements to ensure the existence of upper superposition frame bounds in infinite-dimensional settings, or attempted to characterize the structure of canonical superposition duals in such cases.  Second, establishing additional connections to lapped transform constructions, in particular the tight lapped frames recently proposed in~\cite{Chebira2007}, seem a promising avenue for further investigation.

\appendix

\begin{IEEEproof}[Proof of Lemma~\ref{lem:energyCons}]
To establish~\eqref{eq:energyCons} for all $w \in \mathbb{C}^L$ and $M \geq \operatorname{len}(w)$, expand the left-hand side of~\eqref{eq:energyCons} as
\small
\begin{align}
\sum_{m=0}^{M-1} & \left |\langle x ,\mathcal{M}_{mb}w\rangle \right |^2
= \sum_{m=0}^{M-1} \langle x,\mathcal{M}_{mb}w \rangle \overline{\langle x,\mathcal{M}_{mb} w \rangle}
\notag \\ & = \sum_{t=0}^{L-1} \sum_{t'=0}^{L-1} x[t] \overline{x[t']}\, \overline{w[t]} w[t'] \sum_{m=0}^{M-1} e^{-2\pi i mb(t-t')/L}
\notag \\ & = M \sum_{t=0}^{L-1} \sum_{t'=0}^{L-1} \mathbb{I}_{M\backslash(t-t')}[t-t'] \, x[t] \overline{x[t']}\, \overline{w[t]} w[t']
\label{eq:vanishRoots}
\text{.}
\end{align}
\normalsize
Now consider all $M \geq \operatorname{len}(w)$ that divide $t-t'$; since $\overline{w[t]}w[t'] = 0$ for all $|t-t'|\geq \operatorname{len}(w)$, we need only consider the case $t-t'=0$, whereupon we recover from~\eqref{eq:vanishRoots} the right-hand side of~\eqref{eq:energyCons}.
\end{IEEEproof}

\begin{IEEEproof}[Proof of Lemma~\ref{lem:emw}]
To prove sufficiency, assume that $M_0 = M$, and hence $b_0 = b$.  From~\eqref{eq:energyCons} of Lemma~\ref{lem:energyCons}, it follows that the right-hand side of~\eqref{eq:energyLemma} may be expanded as
\small
\begin{equation}\label{eq:crossTerm}
\!M \!\sum_{t=0}^{L-1} |x[t]|^2 \! \left(|w_p[t]|^2 + |{w_q}'[t]|^2 \!+ 2 \!\operatorname{Re}\{\overline{w_p[t]}{w_q}'[t]\} \right)
\!
\text{,}
\end{equation}
\normalsize
with the rightmost term nonnegative by Definition~\ref{def:superPosition}.  Dropping this term from~\eqref{eq:crossTerm} and applying~\eqref{eq:energyCons} again, this time in the reverse direction, shows that~\eqref{eq:energyLemma} holds for any $M_0 \in \{\max(\operatorname{len}(w_p), \operatorname{len}(w_q)), \ldots, M\}$, thus proving sufficiency.

To prove necessity, assume to the contrary and consider a setting in which $M < \operatorname{len}(w_p + {w_q}') = L$.  Noting that~\eqref{eq:energyLemma} may be stated as $\langle S^{(p,q)} x, x \rangle \leq \langle S^{(p+q)} x, x \rangle$ for positive semi-definite frame operators $S^{(p,q)}$ and $S^{(p+q)}$, assume that elements of the former span $\mathbb{C}^L$ and thus form a frame; hence $\langle S^{(p,q)} x, x \rangle > 0 $ for all nonzero $x \in \mathbb{C}^L$.  However, since $M < L$, the $M$ elements of the latter cannot span $\mathbb{C}^L$.  Hence there exists at least one nonzero $x \in \mathbb{C}^L$ such that $\langle S^{(p+q)} x, x \rangle = 0$, thus contradicting the stated inequality.
\end{IEEEproof}

\begin{IEEEproof}[Proof of Theorem~\ref{thm:testCond}]
To establish the theorem we directly bound the quantity $\sum_{\phi \in \mathscr{F}(I)}  |\langle  x, \phi_{m,n,r} \rangle |^2$ from below.  To begin, observe that
\small
\begin{align*}
&\sum_{\phi \in \mathscr{F}(I)}  |\langle  x, \phi_{m,n,r} \rangle |^2 =
 \sum_{m=0}^{M_{\textrm{g}}-1} \sum_{n,r =0}^{N-1} I[n,r] \left | \langle x,\mathcal{M}_{mb_{\textrm{g}}}\mathcal{T}_{na} w_r \rangle \right |^2
\\ &\!=\!M_{\textrm{g}} \!\!\! \sum_{n,r =0}^{N-1} \!\! I[n,r] \!\! \sum_{t,t'=0}^{L-1} \!x[t] \overline{x[t']} \mathcal{T}_{na}(\overline{w_r[t]}w_r[t']) \mathbb{I}_{M_{\textrm{g}} \backslash (t-t')} [t\!-\!t']
\text{,}
\end{align*}
\normalsize
with the latter expression above obtained by the expansion of~\eqref{eq:vanishRoots}. Now, if $M_{\textrm{g}}$ divides $t-t'$, then $t' = t - k(t)M_{\textrm{g}}$ for some integer $k(t)$, whose domain is deduced by observing that $0 \leq t - k(t)M_{\textrm{g}} \leq L-1$ and $0 \leq t' \leq L-1$ together imply that $k(t) \in \mathscr{K} \triangleq \{ \lceil (t-(L-1))/M_{\textrm{g}}  \rceil, \ldots, \lfloor t/M_{\textrm{g}} \rfloor \}$. Then, a change of variable for $t'$ yields the simplification
\small
\begin{equation*}
    M_{\textrm{g}} \!\!  \sum_{t=0}^{L-1} \sum_{n,r=0}^{N-1} \!\!\!\! \! \sum_{\quad k(t) \in \mathscr{K}} \!\!\! \!\!\!\! x[t] \overline{x[t-kM_{\textrm{g}}]} I[n,r] \mathcal{T}_{na}(\overline{w_r[t]}w_r[t-kM_{\textrm{g}}])
    \text{.}
\end{equation*}
\normalsize
For $k=0$, this quantity can be bounded from below as
\small
\begin{multline*}
\textstyle M_{\textrm{g}} \sum_{t=0}^{L-1} |x[t]|^2 \cdot \sum_{n,r=0}^{N-1} I[n,r]  \left|\mathcal{T}_{na}w_r[t]\right|^2
 \\
 \textstyle \geq  M_{\textrm{g}} \|x\|^2 \cdot \min_{t\in\mathbb{Z}_L} \sum_{n,r=0}^{N-1} I[n,r]  \left|\mathcal{T}_{na}w_r[t]\right|^2 \text{,}
\end{multline*}
\normalsize
with the remaining terms in $\mathscr{K} \setminus \{0\}$ handled as follows.  Invoking the assumption of a real, nonnegative window $w$ to simplify the corresponding expression, observe that the terms $\mathcal{T}_{na}(\overline{w_r[t]}w_r[t-kM_{\textrm{g}}])$ are then everywhere nonnegative.  The sum of the remaining terms can hence be bounded below by
\small
\begin{multline*}
    M_{\textrm{g}} \min_{t \in\mathbb{Z}_L} \sum_{n,r =0}^{N-1} I[n,r] \!\!\!\!  \sum_{k(t) \in \mathscr{K} \setminus \{0\} } \!\!\!\!
    \left(\mathcal{T}_{na}(w_r[t]w_r[t-kM_{\textrm{g}}]) \right)
    \\ \cdot
     \sum_{t''=0}^{L-1} x[t''] \overline{x[t''-kM_{\textrm{g}} ]}
    \\
    \geq - M_{\textrm{g}} \|x\|^2  \min_{t \in\mathbb{Z}_L} \sum_{n,r =0}^{N-1} I[n,r] \!\!\!\! \sum_{k(t) \in \mathscr{K} \setminus \{0\} } \!\!\!\! \mathcal{T}_{na}(w_r[t]w_r[t-kM_{\textrm{g}}]) \text{,}
\end{multline*}
\normalsize
where the second inequality follows by observing that $\sum_{t = 0}^L x[t+s]\overline{x[t]} \geq -\|x\|^2$ for any $f \in \mathbb{C}^L$ and $s \in \mathbb{Z}_L$. Thus, we obtain the claimed results since $\sum_{\phi \in \mathscr{F}(I)}  |\langle  x, \phi_{m,n,r} \rangle |^2$ is bounded from below by $M_{\textrm{g}} \|x\|^2$ times
\small
\begin{equation*}
    \min_{t \in\mathbb{Z}_L} \! \sum_{n,r =0}^{N-1} \!\! I[n,r] \Big( \!\left|\mathcal{T}_{na}w_r[t]\right|^2 - \!\!\!\! \!\!\!\! \sum_{k \in \mathscr{K} \setminus \{0\} } \!\!\!\!\! \mathcal{T}_{na}(w_r[t]w_r[t-kM_{\textrm{g}}]) \Big)
    \text{.}
\end{equation*}
\normalsize
\end{IEEEproof}

\begin{IEEEproof}[Proof of Theorem~\ref{thm:golaStructure}]
To establish the result, first note that $S_{I^{\textrm{g}}}$ is by hypothesis diagonal, and hence by~\eqref{eq:dualWalnutFormula}, we have
\small
\begin{equation}\label{eq:diagFrRef}
\widetilde{\phi}_{m,n,r}[t] =  \frac{\mathcal{M}_{mb_L} \mathcal{T}_{na} w_r[t]}{M_{\textrm{g}} \textstyle{\sum_{n'=0}^{N-1}\sum_{r'=0}^{N-1}} I^{\textrm{g}}[n',r']\left|\mathcal{T}_{n'a}w_{r'}[t]\right|^2}
\text{.}
\end{equation}
\normalsize
If all windows $\{\mathcal{T}_{na}w: n \in \mathbb{Z}_N \}$ have been merged to yield a single superposition window $w_{N-1}$ whose modulates comprise the superposition frame $\mathscr{F}(I^{\textrm{g}})$ of interest, then the constant overlap-add constraint of~\eqref{eq:olaConstraint} applied to~\eqref{eq:diagFrRef} immediately implies the result, as both its numerator and denominator yield constants, whose ratio is in turn $a/(M_{\textrm{g}}\widehat{w}[0])$, with $M_{\textrm{g}} = \max (L, M)$.  Therefore, assume that this is not the case.

To begin, note that~\eqref{eq:diagFrRef}, together with the neighbor-overlap condition of~\eqref{eq:limitOverlap}, implies that $\operatorname{supp} (\widetilde{\phi}_{m,n,r}) \subseteq \mathcal{T}_{na}(L_r \cup C_r \cup R_r)$.  Since these sets are mutually disjoint, we proceed by showing that~\eqref{eq:diagFrRef} agrees with
\small
\begin{equation}\label{eq:lappedEvals}
        \begin{cases}
            \widetilde{\phi}_{m,n,0}[t] & \text{if $t \in \mathcal{T}_{na}L_r$,} \\
            \frac{1}{M_{\textrm{g}}} \frac{a}{\widehat{w}[0]}e^{2\pi i mb_L t / L} & \text{if $t \in \mathcal{T}_{na}C_r$,} \\
            \widetilde{\phi}_{m,n+r+1,0}[t] & \text{if $t \in \mathcal{T}_{na}R_r$,}
          \end{cases}
    \end{equation}
\normalsize
where
\small
\begin{equation*}
    \widetilde{\phi}_{m,n,0}[t] \triangleq \frac{\mathcal{M}_{mb_L} \mathcal{T}_{na}w[t]}{M_{\textrm{g}} \textstyle \sum_{n'=0}^{N-1}\left|w[t-n'a]\right|^2}
    \text{.}
\end{equation*}
\normalsize

We now proceed to show that~\eqref{eq:diagFrRef} evaluates to~\eqref{eq:lappedEvals}.  First, we have that the numerator of~\eqref{eq:diagFrRef} evaluates on $\mathcal{T}_{na}C_r$ to
\small
\begin{align*}
     \mathcal{M}_{mb_L} \mathcal{T}_{na} w_r[t] \, \mathbb{I}_{\mathcal{T}_{na}C_r}[t]
       &= \textstyle \mathcal{M}_{mb_L} \left (\sum_{n'=0}^r \mathcal{T}_{(n+n')a} w[t]
      \right) \, \mathbb{I}_{\mathcal{T}_{na}C_r}[t]
    \\ &= \textstyle  \mathcal{M}_{mb_L}  \left (\sum_{n'=0}^{N-1} \mathcal{T}_{(n+n')a} w[t] \right ) \mathbb{I}_{\mathcal{T}_{na}C_r}[t]
    \\ & =  \textstyle e^{2\pi i mb_L t / L} \left (\frac{\widehat{w}[0]}{a} \right ) \mathbb{I}_{\mathcal{T}_{na}C_r}[t]  \text{,}
\end{align*}
\normalsize
where the second equality follows from the neighbor overlap condition of~\eqref{eq:limitOverlap}, and the third by the overlap-add constraint of~\eqref{eq:olaConstraint} together with the definition of the set $C_r$.  The denominator of~\eqref{eq:diagFrRef} evaluates on this same set $\mathcal{T}_{na}C_r$ to
\small
\begin{align*}
    & \textstyle \left ( M_{\textrm{g}} \textstyle{\sum_{n'=0}^{N-1}\sum_{r'=0}^{N-1}} I^{\textrm{g}}[n',r']\left|\mathcal{T}_{n'a}w_{r'}[t]\right|^2 \right ) \mathbb{I}_{\mathcal{T}_{na}C_r}[t] \\
        & \qquad = \left (M_{\textrm{g}} | \mathcal{T}_{na}w_{r}[t] |^2 \right ) \mathbb{I}_{\mathcal{T}_{na}C_r}[t]  =  \textstyle \left (  M_{\textrm{g}} \frac{\widehat{w}^2[0]}{a^2} \right ) \mathbb{I}_{\mathcal{T}_{na}C_r}[t] \text{,}
\end{align*}
\normalsize
with the first equality following from the fact that no windows other than $\mathcal{T}_{na}w_{r}$ are supported on $\mathcal{T}_{na}C_r$, and the second from~\eqref{eq:olaConstraint}.  Hence we have equality of~\eqref{eq:diagFrRef} and~\eqref{eq:lappedEvals} on $\mathcal{T}_{na}C_r$.

Applying next the neighbor-overlap condition of~\eqref{eq:limitOverlap} and the definition of $L_r$, we observe that the corresponding numerator term of~\eqref{eq:diagFrRef} evaluates to
\small
\begin{align*}
     \mathcal{M}_{mb_L} \mathcal{T}_{na} w_r[t] \, \mathbb{I}_{\mathcal{T}_{na}L_r}[t]
     &=  \mathcal{M}_{mb_L} {\textstyle \left (\sum_{n'=0}^r \mathcal{T}_{(n+n')a} w[t] \right )} \, \mathbb{I}_{\mathcal{T}_{na}L_r}[t]
     \\ & = \mathcal{M}_{mb_L} \mathcal{T}_{na} w[t] \, \mathbb{I}_{\mathcal{T}_{na}L_r}[t] \text{.}
\end{align*}
\normalsize
Evaluating the denominator of~\eqref{eq:diagFrRef} on $\mathcal{T}_{na}L_r$ yields $ (M_{\textrm{g}} \sum_{n'=0}^{N-1} \sum_{r'=0}^{N-1} I^{\textrm{g}}[n',r'] | \mathcal{T}_{n'a} w_{r'}[t] |^2 ) \mathbb{I}_{\mathcal{T}_{na}L_r}[t]$, which may be split into three parts according to index $n$, including the term $\mathcal{T}_{na} w_{r}[t]$ as follows:
\small
\begin{align*}
& M_{\textrm{g}} \textstyle \left( \sum_{r'=0}^{N-1} \sum_{n'=0}^{N-1} I^{\textrm{g}}[n',r'] | \mathcal{T}_{n'a} w_{r'}[t] |^2 \right) \mathbb{I}_{\mathcal{T}_{na}L_r}[t]
\\ & = M_{\textrm{g}} \textstyle \left( \sum_{r'=0}^{N-1} \sum_{n'=0}^{n-1} I^{\textrm{g}}[n',r'] | \mathcal{T}_{n'a} w_{r'}[t] |^2
+ |\mathcal{T}_{na} w_{r}[t] |^2 \right .
\\ & \qquad \qquad + \left . \textstyle \sum_{r'=0}^{N-1} \sum_{n'=n+1}^{N-1} I^{\textrm{g}}[n',r'] | \mathcal{T}_{n'a} w_{r'}[t] |^2 \right) \mathbb{I}_{\mathcal{T}_{na}L_r}[t] \text{.}
\end{align*}
\normalsize
The middle expression of $|\mathcal{T}_{na} w_{r}[t] |^2$ stems from the fact that $I^{\textrm{g}}[n,r]=1$ whenever $\widetilde{\phi}_{\cdot,n,r} \in \widetilde{\mathscr{F}(I^{\textrm{g}})}$, and correspondingly $I^{\textrm{g}}[n,r']=0$ whenever $r' \neq r$, since $I^{\textrm{g}}$ is an admissible selection function.  Moreover, coupled with the assumed neighbor-overlap condition of~\eqref{eq:limitOverlap}, this same property implies that exactly \emph{two} superposition windows are supported on $\mathcal{T}_{na}L_r$, one of which is the superposition window $\mathcal{T}_{na} w_{r}[t]$ isolated in the sum above.

By the superposition construction, it must be the case that the portion of $\mathcal{T}_{na} w_{r}[t]$ supported on $\mathcal{T}_{na}L_r$ takes the form $\mathcal{T}_{na}w[t]$, whereas the portion of the remaining window supported on  $\mathcal{T}_{na}L_r$ takes the form $\mathcal{T}_{(n-1)a} w[t]$.  Thus
\small
\begin{equation}
\label{eq:globalProperty}
\begin{split}
& M_{\textrm{g}} \textstyle \left( \sum_{r'=0}^{N-1} \sum_{n'=0}^{N-1} I^{\textrm{g}}[n',r'] | \mathcal{T}_{n'a} w_{r'}[t] |^2 \right) \mathbb{I}_{\mathcal{T}_{na}L_r}[t]
\\ &  \qquad = M_{\textrm{g}} \textstyle \left ( | \mathcal{T}_{(n-1)a} w[t] |^2  + | \mathcal{T}_{na}w[t] |^2 \right ) \mathbb{I}_{\mathcal{T}_{na}L_r}[t]
\\ &  \qquad = M_{\textrm{g}} \textstyle \sum_{n'=0}^{N-1}|w[t-n'a]|^2 \mathbb{I}_{\mathcal{T}_{na}L_r}[t]
\text{,}
\end{split}
\end{equation}
\normalsize
and~\eqref{eq:diagFrRef} is seen to equal~\eqref{eq:lappedEvals} on the set $\mathcal{T}_{na}L_r$.  The case of $\mathcal{T}_{na}R_r$ proceeds by an identical argument, thereby confirming that~\eqref{eq:diagFrRef} agrees separately on $\mathcal{T}_{na}(L_r,C_r,R_r)$ with the quantities of~\eqref{eq:lappedEvals}, as claimed.

Finally, to complete the proof, observe that the cyclic group setting of $\mathbb{Z}_L$ implies the relation $\widetilde{\phi}_{m,n,0}[t] = \mathcal{M}_{mb_L} \mathcal{T}_{na} \widetilde{\phi}_{0,0,0}[t]$, since for any integer $n$, $\sum_{n'=0}^{N-1}|\mathcal{T}_{n'a}w[t]|^2 \!=\! \sum_{n'=0}^{N-1}|\mathcal{T}_{(n'+n)a}w[t]|^2$.  Applying this relation to~\eqref{eq:lappedEvals}, we obtain the theorem as stated. Note that for \emph{local} $I^l$, the sequence of equalities analogous to those in~\eqref{eq:globalProperty} requires knowledge of the local frequency lattice, which may not be known prior to observing the signal.

\end{IEEEproof}

\begin{IEEEproof}[Proof of Theorem~\ref{thm:dyadicStructure}]
The dyadic superposition frame $\mathscr{F}(I^d)$ represents the set of elements selected from $G_{\cup r}^d(b_\textrm{g})$.  Here we denote its the canonical dual by $\widetilde{\mathscr{F}_{I^d}}(G_{\cup r}^d)$, reflecting the explicit dependence on $G_{\cup r}^d$, and likewise define $\mathscr{F}_{I^d}(\cup \widetilde{G^d_r})$, the set of elements selected by $I^d$ from $\cup_{r \in \mathscr{R}} \widetilde{G^d_r}$.  To establish the result, we must verify the claimed equality
\begin{equation}
    \label{eq:setEquality}
    \widetilde{\mathscr{F}_{I^d}}(G_{\cup r}^d) = \mathscr{F}_{I^d}(\cup \widetilde{G^d_r})
    \text{.}
\end{equation}

We proceed to establish the equality of~\eqref{eq:setEquality} elementwise, noting first that the number of modulates of some $T_{na}w_r$ in $\widetilde{\mathscr{F}_{I^d}}(G_{\cup r}^d)$ is given by $M_{\textrm{g}} = L/b_{\textrm{g}}$ and, by construction, is equal to the number of modulates of the same shifted window in each $G^d_r(b_{\textrm{g}}) = \mathscr{G}(w_r,a,b_{\textrm{g}})$.  We therefore fix $m \in \mathbb{Z}_{M_{\textrm{g}}}$ and $r \in \mathscr{R}$ for the remainder of the proof.

Since the dyadic superposition frame operator $S_{I^d}$ is diagonal, the superposition Walnut formula of~\eqref{eq:dualWalnutFormula} implies that we may write each $\widetilde{\phi}_{m,n,r} \in \widetilde{\mathscr{F}_{I^d}}(G_{\cup r}^d)$ as
\begin{equation}
    \label{eq:leftEquation}
    \widetilde{\phi}_{m,n,r}[t] = \frac{\mathcal{M}_{mb_L}\mathcal{T}_{na}w_r[t]}{M_{\textrm{g}} \textstyle{\sum_{n'=0}^{N-1}\sum_{r'=0}^{N-1}} I^d[n',r'] | \mathcal{T}_{n'a}w_{r'}[t] |^2}
    \text{.}
\end{equation}
Each element $\widehat{\phi}_{m,n'',r} \in \widetilde{G^d_r}$ can be likewise written as:
\begin{equation}
    \label{eq:rightEquation}
    \widehat{\phi}_{m,n'',r}[t] = \frac{\mathcal{M}_{mb_L}\mathcal{T}_{n''a_r}w_r[t]}{\textstyle{M_{\textrm{g}} \sum_{n^*=0}^{N_r-1}} | \mathcal{T}_{n^*a_r}w_{r}[t] |^2 }
    \text{,}
\end{equation}
with $a_r = a(r+1)$, $N_r = L/a_r$, and $0 \leq n'' < N_r$, in accordance with Definition~\ref{defn:dySupFr}.  Note that we ordinarily index modulations of Gabor frame elements by $mb_{\textrm{g}}, m \in \mathbb{Z}_{M_{\textrm{g}}}$, but in~\eqref{eq:rightEquation} we adopt the indexing scheme $mb_L$ for appropriate $m \in \mathbb{Z}_{M_L}$, in order to facilitate its direct comparison to~\eqref{eq:leftEquation}.

In order to establish the equality of sets in~\eqref{eq:setEquality}, we need to show that if $\widetilde{\phi}_{m,n,r} \in \widetilde{\mathscr{F}_{I^d}}(G_{\cup r}^d)$, then the expressions in~\eqref{eq:leftEquation} and~\eqref{eq:rightEquation} are equivalent; i.e.,
\begin{equation}
    \label{eq:elementEquality}
    I^d[m,n,r] = 1 \,\Rightarrow\, \widetilde{\phi}_{m,n,r} = \widehat{\phi}_{m,n,r} \text{.}
\end{equation}

To establish the implication of~\eqref{eq:elementEquality}, we first show that the condition $I^d[m,n,r] = 1$ implies that the numerators of~\eqref{eq:leftEquation} and~\eqref{eq:rightEquation} agree. Since $m$ and $r$ are fixed, this means that for all $n$ such that $I^d[m,n,r] = 1$, there must exist an $0 \leq n'' < N_r$ satisfying
\begin{equation}
    \label{eq:numeratorEquality}
    \mathcal{M}_{mb_L}\mathcal{T}_{na}w_r = \mathcal{M}_{mb_L}\mathcal{T}_{n''a_r}w_r \text{.}
\end{equation}
Equality in~\eqref{eq:numeratorEquality} is achieved when $na = n''a_r = n''(r+1)a$, which clearly holds if $r+1$ divides $n$.  But since $I^d[m,n,r]$ selects elements from $G^d_r = \mathscr{G}(w_r,a(r+1),b_{\textrm{g}})$, then $r+1$ divides $n$ by construction, and~\eqref{eq:numeratorEquality} follows.

The argument for agreement of the denominators is more delicate, because it is \emph{not} true that for all $t \in \mathbb{Z}_L$,
\begin{equation}
    \label{eq:denominatorEquality}
    \textstyle \sum_{n'=0}^{N-1}\sum_{r'=0}^{N-1} I^d[n',r'] | \mathcal{T}_{n'a}w_{r'}[t] |^2  = \sum_{n^*=0}^{N_r - 1} | \mathcal{T}_{n^*a_r}w_{r}[t] |^2
    \text{.}
\end{equation}
Instead, we show that~\eqref{eq:denominatorEquality} holds for all $t \in \operatorname{supp} \left ( \mathcal{T}_{na}w_r \right )$---which, together with~\eqref{eq:numeratorEquality}, is sufficient to establish~\eqref{eq:elementEquality}, and consequently our claimed result.

Let $S_{n,r} \triangleq \operatorname{supp}(\mathcal{T}_{na}w_r)$, with $\mathbb{I}_{S_{n,r}}[t]$ the corresponding indicator function.  Using the same arguments as in the penultimate portion of the proof of Theorem~\ref{thm:golaStructure}, observe that the left-hand side of~\eqref{eq:denominatorEquality} can be decomposed as follows:
\begin{align*}
    & \textstyle \left ( \sum_{r'=0}^{N-1}\sum_{n'=0}^{N-1} I^d[n',r'] |\mathcal{T}_{n'a}w_{r'}[t] |^2 \right ) \mathbb{I}_{S_{n,r}}[t]  \\
    & =  \textstyle  \left(\sum_{r'=0}^{N-1} \sum_{n'=0}^{n-1}  I^d[n',r'] | \mathcal{T}_{n'a}w_{r'}[t]|^2 + \mathcal{T}_{na}w_{r}^2[t] \right. \\
    &  \textstyle  \qquad \qquad \,\left. + \sum_{r'=0}^{N-1} \sum_{n'=n + 1}^{N-1}  I^d[n',r'] | \mathcal{T}_{n'a}w_{r'}[t] |^2 \right) \mathbb{I}_{S_{n,r}}[t] \\
    & = \left ( |\mathcal{T}_{(n-1)a}w[t] |^2 + | \mathcal{T}_{na}w_{r}[t] |^2 + | \mathcal{T}_{(n+r+1)a}w[t] |^2 \right ) \mathbb{I}_{S_{n,r}}[t]
    \text{.}
\end{align*}

Applying the neighbor-overlap requirement of~\eqref{eq:limitOverlap} to the right-hand side of~\eqref{eq:denominatorEquality} then yields
\begin{align*}
    & \textstyle \left ( \sum_{n^*=0}^{N_r - 1} | \mathcal{T}_{n^*a_r}w_{r}[t] |^2 \right ) \mathbb{I}_{S_{n,r}}[t]
    \\ &=  \left ( | \mathcal{T}_{(n-1)a}w[t] |^2 + | \mathcal{T}_{na}w_{r}[t] |^2 + |\mathcal{T}_{(n+r+1)a}w[t]|^2 \right ) \mathbb{I}_{S_{n,r}}[t] \text{,}
\end{align*}
thus establishing the equality of~\eqref{eq:denominatorEquality}, and hence the result.
\end{IEEEproof}

\section*{Acknowledgment}

The authors thank the reviewers for many suggestions that have greatly improved the quality and clarity of this article. In addition, the authors acknowledge helpful discussions co-author with T.~F.~Quatieri of~\cite{Rudoy2008}, as well as very generous long-term dialogues with M.~D\"orfler and L.~Rebollo-Neira that have helped this work to achieve its present form.
\bibliographystyle{IEEEtran}
\bibliography{RudoyBasuWolfe2009b}

\end{document}